\documentclass[11pt]{article}
\usepackage{xcolor}

\usepackage{}
\usepackage{amsfonts}
\usepackage{cite}
\usepackage{mathrsfs,amsfonts,amsmath,amssymb,amsthm,color}
\usepackage{dsfont}
\usepackage{curves}
\usepackage{mathrsfs}
\usepackage{pifont}
\usepackage{graphicx}
\usepackage{float}

\numberwithin{equation}{section}
\numberwithin{figure}{section}
\newtheorem{theorem}{Theorem}[section]
\newtheorem{proposition}{Proposition}[section]
\newtheorem{definition}{Definition}[section]
\newtheorem{lemma}{Lemma}[section]
\newtheorem{corollary}{Corollary}[section]

\renewcommand{\qed}{\hfill\rule{0.5em}{0.809em}}

\def\emptyset{\mbox{{\rm \O}}}

\def\qed{\hfill \rule{4pt}{7pt}}

\begin{document}
\baselineskip 16pt
 %\setcounter{chapter}{1}
 % defining short form------
 \newcommand{\la}{\lambda}
 \newcommand{\si}{\sigma}
 \newcommand{\ol}{1-\lambda}
 \newcommand{\be}{\begin{equation}}
 \newcommand{\ee}{\end{equation}}
 \newcommand{\bea}{\begin{eqnarray}}
 \newcommand{\eea}{\end{eqnarray}}

\baselineskip=0.30in

\vspace*{20mm}

  \begin{center}
 {\Large \bf The $k$-apex trees with minimum augmented Zagreb index}\footnote{Email: liumuhuo@163.com ( M. Liu), pangshumei1979@163.com(Corresponding author, S. Pang),   fbelardo@gmail.com (F. Belardo), akbarali.maths@gmail.com (A. Ali)}

\vspace{10mm}

 {\large \bf Muhuo Liu$^{a}$, Shumei Pang$^{b}$,  Francesco Belardo$^{c}$, Akbar Ali$^{d}$}
\vspace{9mm}

\baselineskip=0.20in
 \baselineskip=0.20in
{ $^{a}$Department of Mathematics and   Research Center for Green\\ Development of Agriculture,  South China Agricultural University,\\
  Guangzhou, 510642, China\\[2mm] $^{b}$International Business College, South China Normal University, \\\small Foshan, 528225, P.R. China\\[2mm]
 $^{c}$Department of Mathematics and Applications ``R. Caccioppoli'',\\ University of Naples ``Federico II'', I-80126 Naples, Italy\\[2mm]
  $^d$Department of Mathematics, College of Science,\\ University of Ha\!'il, Ha\!'il, Saudi Arabia\/}

 \end{center}

\vspace{6mm}

\baselineskip=0.20in

\noindent {\bf Abstract: } {\small } For a connected graph $G$ on at least three vertices, the augmented Zagreb index (AZI) of $G$ is defined as $$AZI(G)=\sum_{uv\in E(G)}\left(\frac{d(u)d(v)}{d(u)+d(v)-2}\right)^{3},$$
being a topological index well-correlated with the formation heat of heptanes and octanes. A $k$-apex tree $G$ is a connected graph admitting a $k$-subset $X\subset V(G)$ such that $G-X$ is a tree, while $G-S$ is not a tree for any $S\subset V(G)$ of cardinality less than $k$. By investigating some structural properties of $k$-apex trees, we identify the graphs minimizing the AZI among all $k$-apex trees on $n$ vertices for $k\ge 4$ and $n\ge 3(k+1)$. The latter solves an open problem posed in [K. Cheng, M. Liu, F. Belardo, {\em Appl. Math. Comput.}, {\bf402} (2021), 126139].
\vspace{2mm}
\begin{flushleft}
\end{flushleft}

\baselineskip=0.25in

\section{Introduction}
Throughout this paper, we restrict ourselves to simple and connected graphs. For a graph $G=(V(G),E(G))$ and $v\in V(G)$, we denote by $N_G(v)$ the set of all neighbors of $v$, and by $|N_G(v)|=d_G(v)$ the degree of $v$ in $G$. A vertex $v$ with $d_G(v)=r$ is said to be an $r$-vertex. As usual, let $\Delta(G)=\max\{d_G(v):v\in V(G)\}$ be the maximum degree of $G$. A 1-vertex is also called a {\bf pendant vertex} (or leaf), a $r$-vertex with $r\geq2$ is called a {\bf non-pendant vertex}, and an $r$-vertex with $r\geq3$ is referred to as a {\bf branching vertex}. An $(s,t)$-edge ($\left(s^{\ge},t^{\ge}\right)$-edge) is an edge whose end vertices have, respectively, degrees $s$ and $t$ (resp., at least $s$ and $t$). We omit the symbol ``$G$'' from notations (as, for example, in $d_G(v)$ or $N_G(v)$) if the graph under consideration is clear from the context.

In Chemical Graph Theory, the chemical and physical properties of a molecule are studied by means of the combinatorial structure of the associated (molecular) graph. In particular, various topological indices based on the structure of a graph, are well-correlated with different chemical and physical characteristics of the corresponding molecule. Among them, the atom-bond connectivity (ABC) index, introduced in  \cite{Gutman} in a slightly different way, is one of such topological indices. The ABC index of a graph $G$ is defined \cite{Est1}
as $$ABC(G)=\sum_{uv\in E(G)} \sqrt{\frac{d(u)+d(v)-2}{d(u)d(v)}}.$$
The ABC index has been associated with the alkanes' heat of formation, and a quantum-chemical explanation for this correlation is provided in \cite{Est1}. From the latter invariant Furtula et al. \cite{Fur2} have derived the generalized ABC index:  $$ABC_{\alpha}(G)=\sum_{uv\in E(G)}\left(\frac{d(u)+d(v)-2}{d(u)d(v)}\right)^{\alpha},$$
where $\alpha$ is a real number and $G$ does not admit any component isomorphic to the 2-vertex path graph whenever $\alpha<0$. The $ABC_{-3}(G)$ is called the augmented Zagreb index (AZI) of $G$. In particular, $$AZI(G)=\sum_{uv\in E(G)}\left(\frac{d(u)d(v)}{d(u)+d(v)-2}\right)^{3}.$$ The authors of \cite{Fur2} observed that the AZI outperforms the ABC index in predicting the heat of formation of heptans and octanes. For further chemical applications of the AZI, ABC index, and several other topological indices, we refer the reader to \cite{BFurtula3,IG1,IG2,y20}.\smallskip

\begin{figure}[h!]
\begin{center} \includegraphics[scale=0.80]{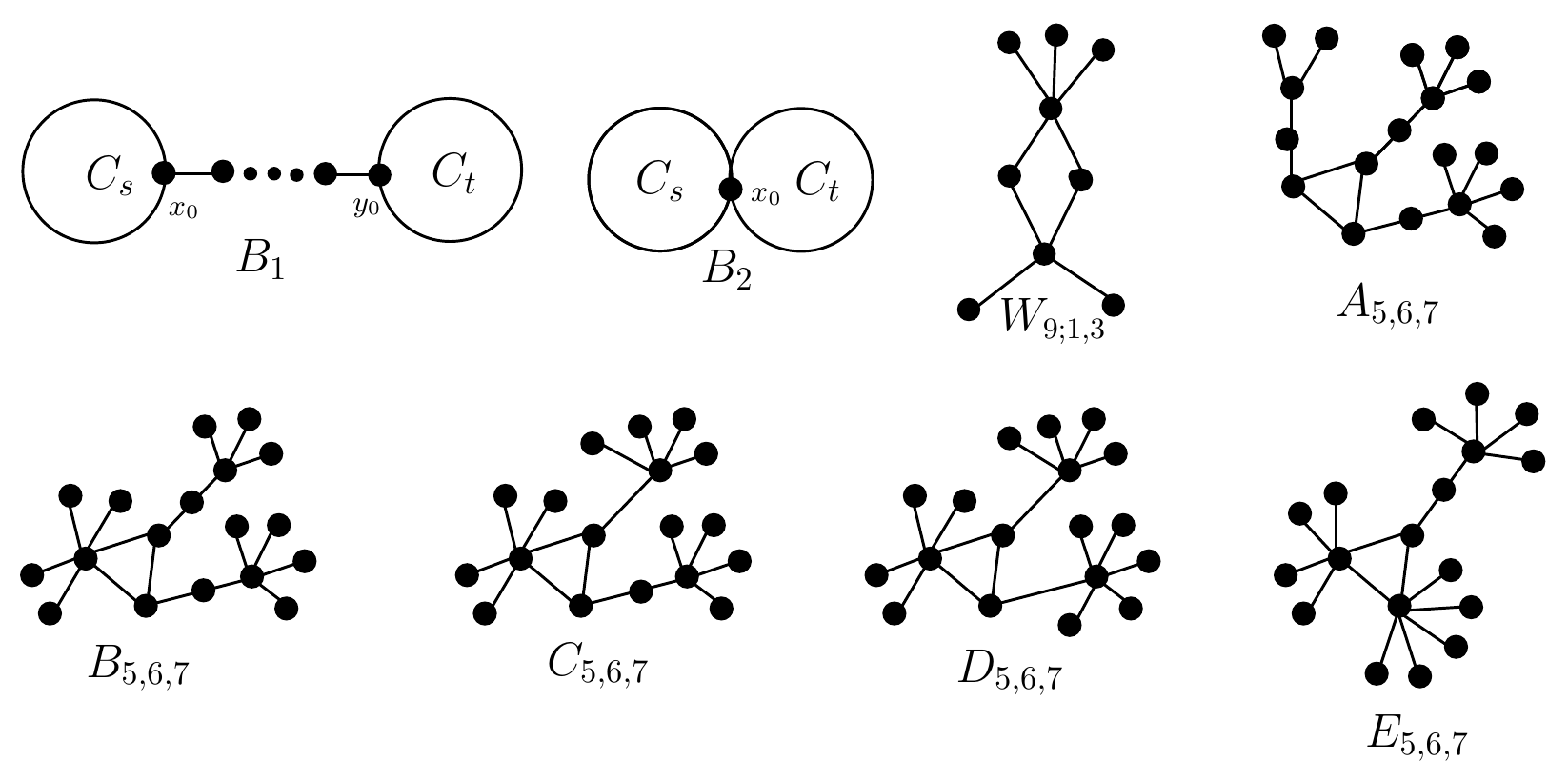}
    \caption{The graphs $B_1$, $B_2$,  $W_{9;1,3}$, $A_{5,6,7}$, $B_{5,6,7}$, $C_{5,6,7}$, $D_{5,6,7}$ and $E_{5,6,7}$.}  \label{21f}\end{center}\end{figure}

We introduce here some additional notation used in this paper.    As usual, $P_n$, $C_n$, and $K_n$ denote the path, the cycle, and the complete graph on $n$ vertices, respectively. Let $K_{s,n-s}$ be the complete bipartite graph on $n$ vertices with $s$ vertices in a partition set and $n-s$ vertices in the other one. Particularly, $K_{1,n-1}$ is the star on $n$ vertices. If $G$ is a  connected graph on $n$ vertices and $n + c - 1$ edges, then $G$ is said to be a $c$-cyclic graph. As customary, for $c = 0, 1, 2, 3$,  $G$ is a tree, unicyclic graph, bicyclic graph, tricyclic graph, respectively. Let  $W_{n;c,s}$ be the $c$-cyclic graph obtained from $K_{2,c+1}$ by attaching $s-1$ and $n-s-c-2$ pendant vertices to the two $(c+1)$-vertices, respectively, where $1\le s\le \frac{1}{2}(n-c-2)$ and $c\ge 1$. For instance, $W_{9;1,3}$ is shown, together with some later described graphs, in Fig. \ref{21f}.
Recently, extremal results about the AZI have been received the attention of several scholars. For example, Furtula et al. \cite{Fur2} showed that the star $K_{1,n-1}$ is the unique graph minimizing the AZI among all trees on $n\geq 3$ vertices. Afterwards, Huang et al. \cite{Huang2} and Zhan et al. \cite{Zhan} independently proved that the graphs obtained from $K_{1,n-1}$ by adding either one edge, or two non-adjacent edges are the graphs minimizing the AZI in the classes of unicyclic and bicyclic graphs on $n\geq5$ vertices, respectively. Ali \cite{Ali1} showed that $W_{3;n,1}$ uniquely minimizes the AZI among all tricyclic graphs of order $n\geq 6$.  Lin et al. \cite{Lin} (see also \cite{Lin3}) characterize the unique tree that maximizes the  AZI in the class of all trees with $n$ vertices for $n\geq 19$, and they also pose the problem of characterizing the extremal  graphs with respect to the AZI among all connected graphs on $n$ vertices and $m$ edges, with $m\geq n$. The latter problem was attacked by  Liu et al. in \cite{MH1}, where it is shown that $W_{n;c,1}$ uniquely minimizes the AZI in the class of all $c$-cyclic graphs on $n$ vertices for $n\geq 2c+1$ and $c\geq 3$. The extremal tree with the maximum and minimum AZI among all trees with $n$ vertices and $k$ pendant vertices were recently identified in \cite{MH22}. For further results on the AZI, the interested reader is referred to the survey paper \cite{Ali2}.

Let $G$ be a connected  graph with at least one cycle. The {\bf base graph} of $G$, denoted by $\mathcal{R}(G)$, is the subgraph of $G$ obtained by recursively deleting the pendant vertices, until there is no pendant vertex left. Evidently, $\mathcal{R}(G)$ is uniquely determined and it keeps the same cyclomatic number of $G$. Conversely, given a graph $G$ with at least a cycle, $G$ is obtained from $\mathcal{R}(G)$ by attaching some hanging trees to the vertices of $\mathcal{R}(G)$. In the forthcoming discussion, by $T_v$ we denote a rooted tree whose root is $v$ (the root is a vertex of $T_v$).

Let $Z_{q}$ be the unique  tree with $q-2$ pendant vertices, one 2-vertex and one $(q-2)$-vertex, where $q\ge 5$, and take $Z_4=P_4$. Hereafter, we refer to the unique pendant vertex   adjacent to the 2-vertex of $Z_{q}$ as the {\bf special pendant vertex} of $Z_{q}$ for $q\ge 5$. For $P_4$,   the two pendant vertices are both  special pendant vertices.

A quasi-tree is a graph $G$ having a vertex $v \in V(G)$ such that $G-v$ is a tree  \cite{Xu1}. Xu et al. \cite{Xu2}  generalized the concept of quasi-tree to $k$-apex trees.

\begin{definition}\label{12d}{\em \cite{Xu2,Xu3}} For $k\geq 1$, a {\bf $k$-apex tree} is a  connected  graph $G$ admitting a $k$-subset $V_k \subset V(G)$, such that $G-V_k$ is a tree,   but $G - X$ is not a tree whenever $X \subset V(G)$ and $|X|< k$. The set $V_k$ is called a {\bf $k$-quasi vertex set} of $G$.
\end{definition}

Denote by $\mathbb{T}^k_n$ the class of $k$-apex trees on $n$ vertices.
In \cite{MH21}, two of the present paper authors and their cooperator determined the unique graphs with the minimum AZI among all graphs of $\mathbb{T}^k_n$ for $1\le k\le3$, and it was observed that providing a complete result for the general case $k\ge 4$ could be a difficult task. In this paper, we give a complete solution to the latter problem, namely, the identification of the $k$-apex trees of order $n$, with $k\geq4$  and $n\ge 3(k+1)$, minimizing the AZI. The main result is obtained by a careful case analysis which stepwise discards $c$-cyclic $k$-apex trees for $c\geq2$, and by identifying the AZI minimizers in the class of unicyclic $k$-apex trees. Hence, let $\mathscr{U}^k_n$ be the class of those unicyclic graphs on $n$ vertices that are  $k$-apex trees.

We now define some graphs in $\mathscr{U}^k_n$, useful to describe the obtained extremal graphs:

\par\noindent
$(i)$   Let $A_{s_1,s_2,s_3}$ be the unicyclic graph with the unique cycle  $C_3:\,=u_1u_2u_3u_1$ such that $T_{u_i}=Z_{s_i}$ and $u_i$ is the special pendant vertex  of $Z_{s_i}$, where $s_i\ge 4$ and   $1\le i\le 3$; \par\noindent
$(ii)$  Let $B_{s_1,s_2,s_3}$ be the unicyclic graph having the unique cycle  $C_3:\,=u_1u_2u_3u_1$ such that $T_{u_1}=K_{1,s_1-1}$ with $u_1$ being the center vertex of $K_{1,s_1-1}$,  and $T_{u_i}=Z_{s_i}$ with  $u_i$ being  the special pendant vertex of $Z_{s_i}$ for   $i\in \{2,3\}$, where $s_1\ge 2$ and $\min\{s_2,s_3\}\ge 4$; \par\noindent
$(iii)$  Let $C_{s_1,s_2,s_3}$ be the unicyclic graph having the unique cycle  $C_3:\,=u_1u_2u_3u_1$ such that $T_{u_1}=K_{1,s_1-1}$ with $u_1$ being the center vertex of $K_{1,s_1-1}$,   $T_{u_2}=K_{1,s_2-1}$ with  $u_2$ being a  pendant vertex of $K_{1,s_2-1}$,  and $T_{u_3}=Z_{s_3}$ with $u_3$ being the special pendant  vertex of $Z_{s_3}$, where $s_1\ge 2$, $s_2\ge 3$ and $s_3\ge 4$.
\par
For the sake of the reader, as an example of the above definitions, we depict the graphs $A_{5,6,7}$, $B_{5,6,7}$ and $C_{5,6,7}$ in Fig. \ref{21f}.
Now, we are ready to state the main result of this paper:

\begin{theorem}\label{1t} Let $G_0$ be a graph with minimum AZI in $\mathbb{T}^k_n$. If $k\ge4$ and $n\ge 3k+3$, then
\begin{align*}G_0=\left\{
\begin{array}{lll}A_{k,\left\lfloor\frac{n-k}{2}\right\rfloor, \left\lceil \frac{n-k}{2}\right\rceil}\,\,\,&\text{if $k\ge 7$},\\[3mm]
B_{6,\left\lceil\frac{n}{2}\right\rceil-3,\left\lfloor\frac{n}{2}\right\rfloor-3}\,\,\,&\text{if $k=6$},\\[3mm]
C_{k,k,n-2k}\,\,\,&\text{if   $4\leq k\leq 5$.}
\end{array}
\right.\end{align*}
\end{theorem}

Here is the remainder of this paper. In Section 2, we provide some bounds for the AZI and some related structural results useful for our investigations. In Section 3, we prove our main result restricted to the class $\mathscr{U}^k_n$ of unicyclic graphs. In Section 4, we discard the $c$-cyclic graphs ($c\geq 4$) in $\mathbb{T}^k_n$ as AZI minimizers. Similarly, we discard bicyclic and tricyclic graphs as AZI minimizers in Sections  5 and 6, respectively. %Finally, in Section 7, by collecting all the results and by discarding $4$-cyclic configurations, we obtain the proof of Theorem \ref{1t}.

\section{Preliminaries}

Let us consider the AZI as a real function on two real variables:$$\Psi(x,y)=\left(\frac{xy}{x+y-2}\right)^{3}.$$
\begin{lemma}\label{21l} Let $x$ and $y$ be two real numbers,  where $\min\{x,y\}\geq 1$. \par\noindent $(i)$  {\em\cite{Huang1}} $\Psi(1,y)$ is decreasing for $y\geq 2$. \par
\noindent $(ii)$ {\em\cite{Huang1}} $\Psi(2,y)=8$ for any real number $y$.\par
\noindent $(iii)$ {\em\cite{Huang1}}  If $y\geq 3$ is fixed, then $\Psi(x,y)>8$ is increasing for $x\geq 3$.\par
\noindent $(iv)$ {\em\cite{MH21}} If $y\ge 4$, then $(y-3)\Psi(1,y-1)>(y-4)\Psi(1,y-2)$.\par
\noindent $(v)$ {\em\cite{MH22}} If $x+y\ge 4$ and $x\ge 2$, then $\Psi(x,y)-\Psi(x-1,y)$ strictly decreases with $x\ge y-1$.  \end{lemma}

From Lemma \ref{21l}$(i)$--$(iii)$, the next result straightly follows.

\begin{corollary}\label{21c}  Let $G$ be a connected graph with $e$ edges, $p$ pendant
vertices,  and maximum degree $\Delta$. If $G$ contains  exactly   $s$ edges such that each one of them has at least one end vertex of degree $2$, then
$$AZI(G)\geq p\left(\frac{\Delta}{\Delta-1}\right)^3+8s+(e-p-s)\left(\frac{9}{4}\right)^3.$$
\end{corollary}

 For a subset $V_1$ of $V(G)$, we denote by $G-V_1$ the graph obtained from $G$ by removing the vertices in $V_1$ and all the incident edges;  if $H$  is an induced subgraph of $G$ with $V(H)\cap V_1=\emptyset$, then $H+V_1$ denotes the subgraph of $G$ induced by  vertex set
$V(H)\cup V_1$. In particular, $G - \{v\}$ is shortened to $G - v$.
 \begin{lemma}\label{22l} Let $S$ be a $k$-quasi vertex set of a $k$-apex tree $G$. If  $X\subset S$ such that $G-X$ is connected, then $G-X$ is a $(k-|X|)$-apex tree.  \end{lemma}
\begin{proof}
 We will prove that $G-X$ is a $(k-|X|)$-apex tree with a $(k-|X|)$-quasi vertex set $S\setminus X$.
By our hypothesis, $G-X-S\setminus X=G-S$ is a tree.
 Contrarily, assume that there exists a subset $Y$ of $V(G-X)$ with $|Y|\le k-|X|-1$
(this inequality is equivalent to $|X\cup Y|\le k-1$ because $X\cap Y=\emptyset$)
such that $G-X-Y$ is a tree. Since $X\cup Y$ is a subset of $V(G)$ and $G-(X\cup Y)=G-X-Y$ is a tree with  $|X\cup Y|\le k-1$, this contradicts the fact that $G$ is a $k$-apex tree.   \end{proof}
\begin{lemma}\label{23l} Let  $u_0v_0$ be an edge of a cycle $C$  pertaining to  a   $k$-apex tree $G$ such that   $u_0$ is a $2$-vertex  of $G$. If   $k\ge 2$, then    $G-u_0$ is a $(k-1)$-apex tree and  $G-u_0v_0$ is  either a $(k-1)$-apex tree or a $k$-apex tree.
\end{lemma}
\begin{proof}
Suppose that $P=v_1v_2\cdots v_{s+1}$ $(s\geq 2)$, is a path joining $v_1$  and $v_{s+1}$  (which need not be distinct) such that $v_1$ and $v_{s+1}$ have degree greater than $2$, while all other vertices $v_2,$ $...,$ $v_s$ are of degree $2$. We may further suppose that   $u_0=v_i$ and $v_0=v_{i+1}$, where $2\le i\le s$. In what follows,  let $C$ be a smallest  cycle of $G$ such that $V(P)\subseteq V(C)$, then  $V(C)$ induces a cycle of $G$.

Clearly,   $G-u_0v_0$ is connected. Let $H_1$ be the component of $G-\{v_1,v_{s+1}\}$ such that $H_1$ contains vertices  $v_2,v_3,\ldots,v_s$. Since $H_1$ is a tree, it holds that
\begin{equation}\label{eq-A-new}
   |V(G)\setminus V(H_1)|\ge k,
\end{equation}
as $G$ is a $k$-apex tree. We prove this lemma by the two below claims.

\par\noindent
{\bf Claim 1.   There exists some $k$-quasi vertex set $S$ of $G$ such that $u_0\in S$.  }
\par\noindent {\bf Proof.}   Let $S$ be a $k$-quasi vertex set of $G$. We may suppose that  $u_0\not\in S$ and so $u_0\in G-S$.

  If $v_1, v_{s+1}\in S$, then
  $G-S$ contains no vertex of $V(G)\setminus V(H_1)$, because $u_0$ is not connected with any vertex of $V(G)\setminus (V(H_1)\cup \{v_1,v_{s+1}\})$ in the graph $G-\{v_1,v_{s+1}\}$. This implies that $V(G)\setminus V(H_1)\subseteq S$, which together with \eqref{eq-A-new} yields $|V(G)\setminus V(H_1)|=k=|S|$.
  Let $H_2=H_1+V(C)\setminus V(H_1)-\{u_0\}$. Then, $|V(H_2)|\ge |V(H_1)|$ with equality if and only if $v_1=v_{s+1}$.
  Since $V(C)$ induces a cycle of $\mathcal {R}(G)$, $H_2$ is also a tree.  Since $G$ is a $k$-apex tree, $|V(G)\setminus V(H_2)|\ge k$. Combining this with  $|V(G)\setminus V(H_2)|\le |V(G)\setminus V(H_1)|\le k$, we have $|V(G)\setminus V(H_2)|=k$. Thus, $S'=V(G)\setminus V(H_2)$ is also a    $k$-quasi vertex set  such that  $u_0\in S'$, and so Claim 1 holds.

  Otherwise, $v_1\not\in S$ or $v_{s+1}\not\in S$. By symmetry, we suppose that $v_1\not\in S$  and so $u_0, v_1\in G-S$. It is easily to see that $S''=S\setminus \big(V(H_1)\cup V(C)\big)\cup \{u_0\}$ is also a vertex set of $V(G)$ such that $G-S''$ is also a tree. Since  $G-S$ is a tree and $C$ is a cycle of $G$, $S$ must  contain  at least one vertex of $C$ and thus   $|S''|\le k$. Combining this with  $G$ being a $k$-apex tree, we can conclude that $k\le |S''|\le k$,  as $G-S''$ is a tree. Therefore,   $S''$ is a $k$-quasi vertex set of $G$ with $u_0\in S''$, and so Claim 1 also holds.    \qed

 By Lemma \ref{22l} and Claim 1, we can conclude that $G-u_0$ is a $(k-1)$-apex tree. Suppose that $N_G(u_0)=\{v_0,w_0\}$.

\par\medskip \noindent
{\bf Claim 2.      $G-u_0v_0$ is   either a  $(k-1)$-apex tree or a $k$-apex tree.}
\par\noindent {\bf Proof.} If  $G-u_0v_0$ contains a vertex set $X$ with $|X|\le k-2$ such that $H_2=G-u_0v_0-X$ is a tree, then $X_1=X\cup \{u_0\}$ is a subset of $V(G)$ such that $G-X_1$ is a tree when $u_0, v_0\in V(H_2)$ and $X$ is a subset of $V(G)$ such that $G-X$ is a tree when $u_0\not\in V(H_2)$ or $v_0\not\in V(H_2)$, which is contrary with the fact that $G$ is a $k$-apex tree.

By Claim 1, let $Y$ be a $(k-1)$-quasi vertex set of $G-u_0$. Then,  $G-u_0v_0-Y-u_0=G-Y-u_0$ is a  tree. Since $|Y\cup \{u_0\}|=k$, we can conclude that we can obtain a tree from $G-u_0v_0$ by deleting at most $k$ vertex. Furthermore, from the former arguments, any deletion of at most $k-2$ vertices from $G-u_0v_0$ cannot be a tree. Thus,    $G-u_0v_0$ is  either a $(k-1)$-apex tree or a $k$-apex tree. This completes the proof of Claim 2.  \end{proof}

\begin{lemma}\label{24l} Let $u_0$ be a $2$-vertex  of the  connected graph  $G$. If $u_0$ is a vertex of some cycle $C$  of $G$, there exists an edge $u_0v_0$ of $C$ such that  $AZI(G)-AZI(G-u_0v_0)\ge 8.$
\end{lemma}
\begin{proof}  Suppose that $C$ is a  cycle of $G$ such that  $u_0v_0\in E(C)$,  $w_0u_0\in E(C)$ and  $vv_0\in E(C)$, where $w_0=v$ is also permitted. Since $w_0, v\in V(C)$, it holds that $\min\{d_G(w_0),d_G(v)\}\ge2$.

\par\noindent{\bf Case 1. $d_G(v_0)=2$.}
   By Lemma \ref{21l}(i),   \begin{align}\label{21e}AZI(G)-AZI(G-u_0v_0)= 24-\left(\frac{d_G(w_0)}{d_G(w_0)-1}\right)^3-\left(\frac{d_G(v)}{d_G(v)-1}\right)^3\ge 24-16=8.\end{align}

\par\noindent{\bf Case 2. $d_G(v_0)\ge 3$.}

By Lemma \ref{21l}(i) and Lemma \ref{21l}(iii),  for $d_G(z)\ge3$ we have
\begin{align*}
\Psi\big(d_G(v_0),d_G(z)\big)- \Psi\big(d_G(v_0)-1,d_G(z)\big)&>0=\Psi\big(d_G(v_0),2\big)-\Psi\big(d_G(v_0)-1,2\big)\\
&>\Psi\big(d_G(v_0),1\big)-\Psi\big(d_G(v_0)-1,1\big)
\end{align*}
 and hence
\begin{eqnarray}\label{22e}
% \nonumber to remove numbering (before each equation)
 \nonumber &&AZI(G)-AZI(G-u_0v_0)\\  \nonumber &\ge&  16-\Psi\big(d_G(w_0),1\big)+\sum_{z\in N_G(v_0)\setminus \{u_0,v\}}\Big[\Psi\big(d_G(v_0),d_G(z)\big) -\Psi\big(d_G(v_0)-1,d_G(z)\big)\Big]\\
  &\ge& 16-\Psi\big(d_G(w_0),1\big)+\big(d_G(v_0)-2\big)\Big[\Psi\big(d_G(v_0),1\big) -\Psi\big(d_G(v_0)-1,1\big)\Big].
\end{eqnarray}

%  \begin{align}\nonumber&AZI(G)-AZI(G-u_0v_0)\\\nonumber\ge & 16-\Psi\big(d_G(w_0),1\big)+\sum_{z\in N_G(v_0)\setminus \{u_0,v\}}\Big[\Psi\big(d_G(v_0),d_G(z)\big)-\big(d_G(v_0)-1,d_G(z)\big)\Big]\\\ge & 16-\Psi\big(d_G(w_0),1\big)+\big(d_G(v_0)-2\big)\Big[\Psi\big(d_G(v_0),1\big)-\big(d_G(v_0)-1,1\big)\Big].\end{align}

Let \begin{align*}f(x)=\big(x-2\big)\Big[\Psi\big(x,1\big)-\Psi\big(x-1,1\big)\Big]=\big(x-2\big)\left[\left(\frac{x}{x-1}\right)^3-\left(\frac{x-1}{x-2}\right)^3\right],~\text{where $x\ge 3$.}\end{align*}\label{23e}
Since  $x\ge 3$, we have
$(x-1)^4(x-2)^3f'(x)=3x^5 - 12x^4 + 9x^3 + 18x^2 - 25x + 4>0$, which implies that the function $f$ is a strictly increasing function for $x\ge 3$. By Case 1, we may suppose that $d_G(w_0)\ge 3$.  Combining this with Lemma \ref{21l}(i) and \eqref{22e}, we have
\begin{align*}AZI(G)-AZI(G-u_0v_0)\ge  16-\left(\frac{3}{2}\right)^3+\left(\frac{3}{2}\right)^3-8=8.\end{align*}
This completes the proof.
\end{proof}

Let $u_0$ be a vertex in the base graph of  $G$.  If $|V(T_{u_0})|=q\geq 3$, then denote by $G^*(u_0)$ the graph obtained from $G$ by deleting all vertices of $V(T_{u_0})\setminus \{u_0\}$ and then   adding one new edge between $u_0$ and a pendant vertex of $K_{1,q-2}$.   From the definition, it is easily checked that \begin{align*}\text{If $G$ is a $k$-apex tree, then $G^*(u_0)$ is also a $k$-apex tree for $k\geq 2$.}\end{align*}
\begin{proposition}\label{21p}{\em \cite{MH21}}
Let
 $u_0\in V( \mathcal{R}(G))$  with  $N_G(u_0)\cap V(\mathcal{R}(G))=\{v_0,w_0\}$
and $d_G(v_0)\geq d_G(w_0)\geq 3.$   Suppose that  $|V(T_{u_0})|=q\geq 3$    and $V(T_{u_0})\setminus \{u_0\}$ contains exactly $t$ non-pendant vertices. If $G\neq G^*(u_0)$, then we have $AZI(G)> AZI\big(G^*(u_0)\big)$ if one of the following conditions holds:
\begin{itemize}
{\baselineskip=0.08in
\item[$(i)$] $t=0$, $d_G(v_0)=3$ and $q\geq 7$;\,\,
 $(ii)$  $t=0$, $d_G(v_0)\geq 4$ and $q\geq 5$;
\item[$(iii)$] $t=1$,    $d_G(u_0)=3$ and $q\geq 6$;\,\,
$(iv)$  $t=1$, $d_G(u_0)\geq 4$ and $q\geq 5$;\,\,\,\,$(v)$  $t\geq 2$.}
\end{itemize}
\end{proposition}
At the end of this section, we shall list the formulas  of AZI for different graphs, which will be used  in the sequel  sections.  \begin{align}\label{41e}&AZI(W_{n;c,s})=(n-c-2-s)\left(\frac{n-s-1}{n-2-s}\right)^3+(s-1)\left(\frac{s+c}{s+c-1}\right)^3+16(c+1),
\\[3mm]\label{42e}&AZI(A_{k,k,n-2k})=\frac{5259}{64}+\frac{(n-2k-2)^3}{(n-2k-3)^2}+2\frac{(k-2)^3}{(k-3)^2},
\\[3mm]\label{43e}&AZI(C_{4,4,n-8})=\frac{5289}{64}+\frac{(n-10)^3}{(n-11)^2},\,\,AZI(C_{5,5,n-10})=\frac{6961}{78}+\frac{(n-12)^3}{(n-13)^2},
\\[3mm]\label{44e}&AZI\left(A_{k,\left\lfloor\frac{n-k}{2}\right\rfloor, \left\lceil \frac{n-k}{2}\right\rceil}\right)=\frac{5259}{64}+\frac{(k-2)^3}{(k-3)^2}+\left\{
\begin{array}{ll}
\frac{(n-k-4)^3}{(n-k-6)^2} \hspace*{8pt} \text{if $n-k$ is even},\\[4mm]
\frac{(n-k-3)^3}{2(n-k-5)^2}+\frac{(n-k-5)^3}{2(n-k-7)^2} \hspace*{8pt} \text{if $n-k$ is odd},\\
\end{array}
\right.\end{align}\begin{align}\label{45e} &AZI\left(B_{6,\left\lceil\frac{n}{2}\right\rceil-3,\left\lfloor\frac{n}{2}\right\rfloor-3}\right)
=\frac{7088}{81}+\left\{
\begin{array}{ll}
\frac{(n-10)^3}{(n-12)^2} \hspace*{8pt} \text{if $n$ is even},\\[4mm]
\frac{(n-9)^3}{2(n-11)^2}+\frac{(n-11)^3}{2(n-13)^2} \hspace*{8pt} \text{if $n$ is odd}.\\
\end{array}
\right.\end{align}

 \section{The case of unicyclic graphs}
In this section, we consider $k$-apex trees in $\mathscr{U}_n^k$, and we prove Theorem \ref{1t} under the latter restriction. The main result of this section is

\begin{proposition}\label{31p}Let $G_0$ have the minimum AZI in $\mathscr{U}_n^k$. If $n\ge 3k+3$, then
\begin{align*}G_0=\left\{
\begin{array}{lll}A_{k,\left\lfloor\frac{n-k}{2}\right\rfloor, \left\lceil \frac{n-k}{2}\right\rceil}\,\,\,&\text{if $k\ge 7$},\\[3mm]
B_{6,\left\lceil\frac{n}{2}\right\rceil-3,\left\lfloor\frac{n}{2}\right\rfloor-3}\,\,\,&\text{if $k=6$},\\[3mm]
C_{k,k,n-2k}\,\,\,&\text{if   $4\leq k\leq 5$.}
\end{array}
\right.\end{align*}
 \end{proposition}
To show Proposition \ref{31p}, we need to introduce a series of lemmas.
\begin{lemma}\label{31l}{\em \cite{MH21}}
Let $G\in \mathscr{U}_n^k$, where $k\geq 2$. If $C_g$ is the unique cycle of  $G$, then $\min\big\{|V(T_u)|:\,\,u\in V(C_g)\big\}=k$ and   $G[V_k]$ is a root tree for each $k$-quasi vertex set $V_k$ of $G$.
 \end{lemma}

Let $G^*_{s_1,s_2}$ be the   graph obtained from a connected graph $G$  with $w_1, w_2\in V(G)$ and two stars $K_{1,s_1-2}$ and $K_{1,s_2-2}$ by adding
 an edge between $w_i$ and one pendant vertex   of $K_{1,s_i-2}$ for each $i\in \{1,2\}$, where  $w_1=w_2$ is also permitted.
\begin{lemma}\label{32l} If $s_1\ge s_2+2$ and $s_2\ge 4$, then $AZI(G^*_{s_1,s_2})>AZI(G^*_{s_1-1,s_2+1}).$\end{lemma}
 \begin{proof} By the definition of $G^{*}_{s_1,s_2}$, we have
  \begin{align*}AZI(G^{*}_{s_1,s_2})-AZI(G^{*}_{s_1-1,s_2+1})&=\frac{(s_1-2)^3}{(s_1-3)^2}+\frac{(s_2-2)^3}{(s_2-3)^2}-\frac{(s_1-3)^3}{(s_1-4)^2}
 -\frac{(s_2-1)^3}{(s_2-2)^2}\\&=\frac{(s_1-s_2-1)f(s_1)}{(s_1-3)^2(s_1-4)^2 (s_2-2)^2 (s_2-3)^2}, \end{align*}
where $f(s_1)=(3s^2_2 - 13s_2 + 13)s^3_1 +(s_2 - 13)(3s^2_2- 13s_2 + 13)s^2_1-(19s^3_2 - 238s^2_2 + 756s_2 - 672)s_1 + 29s^3_2 - 319s^2_2 + 960s_2 - 828$. Since $s_1\ge s_2+2$ and $s_2\ge 4$, we have $f'''(s_1)=6(3s^2_2 - 13s_2 + 13)>0$ and thus $f''(s_1)\ge f''(s_2+2)=2(4s_2 - 7)(3s^2_2 - 13s_2 + 13)>0$. Therefore, $f'(s_1)\ge f'(s_2+2)=(s_2 - 1)(s_2 - 2)(15s^2_2 - 69s_2 + 76)>0$, which implies that $f(s_1)\ge f(s_2+2)=2(s_2-2)^2(3s^3_2 - 15s^2_2 + 21s_2 - 7)>0$. \end{proof}
\begin{corollary}\label{31c}If $n\ge 3k+1$ and $k\ge 5$, then $$ AZI\left(A_{k,\left\lfloor\frac{n-k}{2}\right\rfloor, \left\lceil \frac{n-k}{2}\right\rceil}\right)<AZI\left(A_{k-1,\left\lfloor\frac{n-k+1}{2}\right\rfloor, \left\lceil \frac{n-k+1}{2}\right\rceil}\right)< AZI(A_{k-1,k-1,n+2-2k}).$$   \end{corollary}
\begin{proof} Since $k-1\ge 4$ and $n+2-2k>k-1+2$, the second inequality follows from Lemma \ref{32l}.
Since $\left\lceil \frac{n-k+1}{2}\right\rceil\ge \frac{n-k+1}{2}\ge k-1+2$, we have $$AZI\left(A_{k-1,\left\lfloor\frac{n-k+1}{2}\right\rfloor, \left\lceil \frac{n-k+1}{2}\right\rceil}\right)>AZI\left(A_{k,\left\lfloor\frac{n-k+1}{2}\right\rfloor, \left\lceil \frac{n-k-1}{2}\right\rceil}\right)=AZI\left(A_{k,\left\lfloor\frac{n-k}{2}\right\rfloor, \left\lceil \frac{n-k}{2}\right\rceil}\right)$$ by Lemma \ref{32l}.    \end{proof}

\begin{lemma}\label{33l} Let  $u\in V(\mathcal {R}(G))$  such that $T_{u}=Z_{q}$ with $u$ being the special pendant vertex of $Z_{q}$ and $G^*$ be the graph obtained from  $G$  by replacing  $Z_{q}$ with    $K_{1,q-1}$ such that $u$ is a pendant vertex of $K_{1,q-1}$.
If  $4\le q\le 5$ and $d_G(u)=3$, then    $AZI(G)>AZI(G^*),$ where $G\in   \mathscr{U}^k_n$ implies that $G^*\in   \mathscr{U}^k_n$.\end{lemma}
 \begin{proof} Clearly, if $G\in   \mathscr{U}^k_n$, then  $G^*\in   \mathscr{U}^k_n$. Since $d_G(u)=d_{G^*}(u)=3$ and $4\le q\le 5$, it is easily checked that    \begin{align*}AZI(G)-AZI(G^*)\ge 16+2\Psi(3,1)-3\Psi(4,1)-\Psi(3,4)>1.8,\end{align*}    thus the result holds.  \end{proof}

\begin{lemma}\label{34l} If $s_3\ge 5$ and $s_1\ge 3$, then    $AZI(C_{s_1,s_2,s_3})>AZI(C_{s_1-1,s_2,s_3+1})$.  \end{lemma}
\begin{proof} By Lemma \ref{21l}, it follows that
\begin{eqnarray*}
% \nonumber to remove numbering (before each equation)
  AZI(C_{s_1,s_2,s_3})-AZI(C_{s_1-1,s_2,s_3+1}) &=& 2\Psi(3,s_1+1)-2\Psi(3,s_1)+\big((s_1-1)\Psi(1,s_1+1)\\
   &&-(s_1-2)\Psi(1,s_1)\big) +\big((s_3-3)\Psi(1,s_3-2) \\
   && -(s_3-2)\Psi(1,s_3-1)\big)\\
   &>&  2\Psi(s_1+1,3)-2\Psi(s_1,3)-\Psi(1,s_3-2) \\
   &\ge& 2\Psi(3,4)-2\Psi(3,3)-\Psi(1,3)>0.
\end{eqnarray*}
This completes the proof.   \end{proof}

In the rest of this section, {\bf  we always suppose that  $G_0$ has the minimum AZI in $\mathscr{U}^k_n$ and  $C_g$ is the unique cycle of $G_0$}, where $k\ge 2$.
\begin{lemma}\label{35l} Let  $u_1, u_2, u_3$ be three vertices of $C_g$ such that $u_2\in N_{G_0}(u_1)\cap N_{G_0}(u_3)$   and let $G^*$ be the unicyclic graph obtained from $G_0$  by  replacing $T_{u_2}$ with $K_{1,q-1}$ such that $u_2$ is the center vertex of $K_{1,q-1}$, where $3\le q=|V(T_{u_2})|\le 6$.   If $d_{G_0}(u_1)=d_{G_0}(u_3)=3$, then $G_0=G^*$.\end{lemma}
   \begin{proof} By contradiction, we assume that $G_0\ne G^*$. Suppose that $V(T_{u_2})\setminus \{u_2\}$ contains exactly $t$ non-pendant vertices.

   If $T_{u_2}=Z_q$ with $u_2$ being the special pendant vertex of $Z_q$, then $q=6$ by the choice of $G_0$ and Lemma \ref{33l}. Since $d_{G_0}(u_1)=d_{G_0}(u_3)=3$, we have
       \begin{align*}AZI(G_0)-AZI(G^*)=16+3\Psi(4,1)+2\Psi(3,3)-5\Psi(7,1)-2\Psi(3,7)>1.75,\end{align*}    contrary with the choice of $G_0$, as $G^*\in \mathscr{U}^k_n$.

    If $T_{u_2}=K_{1,q-1}$ with $u_2$ being a pendant vertex of $K_{1,q-1}$, then
       \begin{align*}AZI(G_0)-AZI(G^*)&\ge 4\Psi(5,1)+\Psi(3,5)+2\Psi(3,3)-5\Psi(7,1)-2\Psi(3,7)>2,\end{align*}    for $3\le q\le 6$. This is contrary with the choice of $G_0$.

     If  $T_{u_2}\ne  Z_q$ and $T_{u_2}\ne K_{1,q-1}$ with $u_2$ being a pendant vertex of $K_{1,q-1}$, then  Proposition \ref{21p} and the choice of $G_0$ imply that $T_{u_2}=P_4$ with $u_2$ being a 2-vertex of $P_4$. In this case,  \begin{align*}AZI(G_0)-AZI(G^*)= 16+2\Psi(4,3)+\Psi(4,1)-2\Psi(5,3)-3\Psi(5,1)>8.5,\end{align*}   contrary with the choice of $G_0$.
                  \end{proof}

 \begin{lemma}\label{36l}  If   $g\ge 4$, then there exists  $G^*\in \mathscr{U}^k_n$  such that the unique cycle of $G^*$ is  $C_{g-1}$ and  $AZI(G_0)=AZI(G^*)$.      \end{lemma}
\begin{proof}   Let $S=\{w_0,w_1,\ldots,w_{k-1}\}$ be a $k$-quasi vertex set of $G_0$. By   Lemma \ref{31l} and $k\ge 2$, we may suppose that  $G_0[S]$ is a root tree with root $w_0$ and each vertex of $C_g$ is a branching vertex.
   Let $y_1, y_2, y_3,y_4$ be four vertices of $C_g$ such that $y_1\in N_{G_0}(w_0)\cap N_{G_0}(y_2)$, $y_2\in N_{G_0}(y_1)\cap N_{G_0}(y_3)$ and $y_3\in N_{G_0}(y_2)\cap N_{G_0}(y_4)$, where $y_4=w_0$ implies that $g=4$. We choose  $z_i\in N_{G_0}(y_i)\setminus V(C_g)$ such that $d_{G_0}(z_i)$ is as small as possible, where $1\le i\le 4$. Let $G_1=G_0+w_0y_2+y_1z_2-y_1w_0-y_2z_2$. Since  $g\ge 4$, then  $G_1$ is also a $k$-apex tree with girth $g-1$ by Lemma \ref{31l}.\smallskip

\par\noindent{\bf Case 1.  $d_{G_0}(y_2)\le d_{G_0}(y_1)$.}

If $d_{G_0}(z_2)\le 2$, then Lemma \ref{21l} implies that
 \begin{eqnarray*}
% \nonumber to remove numbering (before each equation)
  AZI(G_0)-AZI(G_1) &=& \Psi\big(d_{G_0}(w_0),d_{G_0}(y_1)\big)+\Psi\big(d_{G_0}(z_2),d_{G_0}(y_2)\big) -\Psi\big(d_{G_0}(w_0),d_{G_0}(y_2)\big) \\
   &&-\Psi\big(d_{G_0}(z_2),d_{G_0}(y_1)\big) \\
   &\ge& \Psi\big(d_{G_0}(w_0),d_{G_0}(y_1)\big)- \Psi\big(d_{G_0}(w_0),d_{G_0}(y_2)\big)\ge 0,
\end{eqnarray*}
as  $y_1, y_2, w_0$ are three branching vertices of $G_0$. By the choice of   $G_0$, we have \begin{align}\label{31e}\text{  $d_{G_0}(y_1)=d_{G_0}(y_2)$ and $AZI(G_0)=AZI(G_1)$.}\end{align}

If $d_{G_0}(z_2)\ge 3$, then by the choice of $z_2$, we can conclude that $d_{G_0}(w)\ge 3$ for each vertex $w\in N_{G_0}(y_2)\setminus V(C_g)$.
 Next, we  suppose that  $d_{G_0}(y_2)<d_{G_0}(y_1),$
as $AZI(G_0)=AZI(G_1)$ for $d_{G_0}(y_2)=d_{G_0}(y_1)$.
 By the choice of $G_0$, Proposition \ref{21p} and the choice of $z_2$, we can conclude that \begin{align}\label{32e}\text{ $N_{G_0}(y_2)\setminus V(C_g)=\{z_2\}$ and   $d_{G_0}(y_2)=3<4\le q\le 5$, where  $T_{y_2}=K_{1,q-1}$.}\end{align} Let $d_{G_0}(w_0)=x$ and $d_{G_0}(y_1)=y$. Since $3=d_{G_0}(y_2)<d_{G_0}(y_1),$  we have  $y\ge 4$ and $x\ge 3$. \smallskip

  \par\noindent{\bf Subcase 1.1.     $3\le x\le 4$.}
 Since $d_{G_0}(y_1)=y\ge 4$ and $AZI(G_0)$ is minimum, Proposition \ref{21p} implies that  either  $T_{y_1}=K_{1,r-1}$ with $3\le r\le 6$ and $y_1$ being the center vertex of $T_{y_1}$  or $T_{y_1}=P_4$ with $y_1$ being a 2-vertex of $P_4$. Let $G_2=G_0+w_0y_2+y_1z-y_1w_0-y_1y_2$, where $z$ is a pendant vertex of $T_{y_2}$. Then, $G_2$ is also a $k$-apex tree with girth $g-1$ by Lemma \ref{31l}.

 We only prove the case of $x=3$, as the case $x=4$ can be shown similarly.

We first suppose that   $T_{y_1}=K_{1,r-1}$ with $3\le r\le 6$. By \eqref{32e} and $4\le q\le 5$, we have
\begin{eqnarray*}
% \nonumber to remove numbering (before each equation)
  AZI(G_0)-AZI(G_2) &=& 2\Psi\big(r+1,3\big)+(r-1)\Psi\big(r+1,1\big)+(q-2)\Psi\big(q-1,1\big) \\
   && -(r-1)\Psi\big(r,1\big)-16-(q-3)\Psi\big(q-1,1\big)-\Psi(3,3) \\
   &\ge& 2\Psi\big(r+1,3\big)+(r-1)\Psi\big(r+1,1\big)-(r-1)\Psi\big(r,1\big)+\Psi\big(4,1\big)\\
   && -16-\Psi(3,3)\\
   &=& 2\left(\frac{3(r+1)}{r+2}\right)^3+(r-1)\left(\frac{r+1}{r}\right)^3-(r-1)\left(\frac{r}{r-1}\right)^3+\left(\frac{4}{3}\right)^3\\
   && -\left(\frac{9}{4}\right)^3-16\\
   &>& 2\left(\frac{3(r+1)}{r+2}\right)^3+(r-1)\left(\frac{r+1}{r}\right)^3-(r-1)\left(\frac{r}{r-1}\right)^3-\frac{51}{2}\\
   &=& \frac{g(r)}{2\,(r-1)^2\, r^3\,(r+2)^3}>0,
\end{eqnarray*}

%\begin{align*}&AZI(G_0)-AZI(G_2)\\=&2\Psi\big(r+1,3\big)+(r-1)\Psi\big(r+1,1\big)+(q-2)\Psi\big(q-1,1\big)\\&-
%  (r-1)\Psi\big(r,1\big)-16-(q-3)\Psi\big(q-1,1\big)-\Psi(3,3)
%  \\\ge &2\Psi\big(r+1,3\big)+(r-1)\Psi\big(r+1,1\big)-(r-1)\Psi\big(r,1\big)+\Psi\big(4,1\big)
% -16-\Psi(3,3)
%\\= &2\left(\frac{3(r+1)}{r+2}\right)^3+(r-1)\left(\frac{r+1}{r}\right)^3-(r-1)\left(\frac{r}{r-1}\right)^3+\left(\frac{4}{3}\right)^3
%-\left(\frac{9}{4}\right)^3-16\\> &2\left(\frac{3(r+1)}{r+2}\right)^3+(r-1)\left(\frac{r+1}{r}\right)^3-(r-1)\left(\frac{r}{r-1}\right)^3-25.5\\=&  \frac{g(r)}{2\,(r-1)^2\, r^3\,(r+2)^3}>0,\end{align*}
where $f(r)=57r^8 - 102r^7 - 303r^6 + 228r^5 + 300r^4 - 230r^3 + 36r^2 - 24r - 16.$ This is contrary with the choice of $G_0$.

We second suppose that   $T_{y_1}=P_4$ with  $y_1$ being a 2-vertex of $P_4$. Then,  $d_{G_0}(y_1)=4$ and  \begin{align*}AZI(G_0)-AZI(G_2)&=2\Psi(3,4)+\Psi(4,1)+\Psi(q-1,1)-\Psi(3,3)-16-\Psi(3,1)\\&\ge 2\Psi(3,4)+2\Psi(4,1)-\Psi(3,3)-16-\Psi(3,1)>1.5\,\,\,\,  \text{by Lemma \ref{21l}},\end{align*}
contrary with the choice of $G_0$.

 \par\noindent{\bf Subcase 1.2.    $x\ge 5$.} By \eqref{32e}, we have $4\le q\le 5$.   Since $y\ge 4$, we have  \begin{align*}\Psi\big(3,3\big)
-\Psi\big(y,3\big)-\Psi\big(3,4\big)+\Psi\big(y,4\big)>\Psi\big(y,4\big)
-\Psi\big(y,3\big)-2.5=\frac{ f(y)}{2(y+2)^3(y+1)^3}>0,\end{align*}
where $f(y)=69y^6 + 15y^5 - 429y^4 - 619y^3 - 330y^2 - 180y - 40$ and thus
\begin{align}\nonumber AZI(G_0)-AZI(G_1)&=\Psi\big(x,y\big)+\Psi\big(3,q-1\big)-
  \Psi\big(x,3\big)-\Psi\big(y,q-1\big)\\\label{liun1}&\ge \Psi\big(x,y\big)+\Psi\big(3,4\big)-
  \Psi\big(x,3\big)-\Psi\big(y,4\big).\end{align}
   It is not difficult to see that $\Psi\big(x,y\big)-\Psi\big(x,3\big)$ is strictly increasing with $x\ge 5$ and so
  \begin{align*} AZI(G_0)-AZI(G_1)   \ge \Psi\big(5,y\big)+\Psi\big(3,4\big)-
  \Psi\big(5,3\big)-\Psi\big(y,4\big)=\frac{(y-3)g(y)}{1000(y+2)^3(y+3)^3}\end{align*} by \eqref{liun1},
  where $g(y)=59199y^5 + 324582y^4 + 578253y^3 + 457454y^2 + 367404y + 129672>0$.   This is  contrary with the choice of $G_0$.

\par\noindent{\bf Case 2.  $d_{G_0}(y_2)>d_{G_0}(y_1)$.} Let $G_3=G_0+y_1y_3+y_2z_1-y_2y_3-y_1z_1$, where $z_1$ has minimum degree among all vertices of
 $N_{G_0}(y_1)\setminus V(C_g)$.  Since  $g\ge 4$, $G_3$ is also a $k$-apex tree with girth $g-1$ by Lemma \ref{31l}, as $w_0\not\in \big\{y_1,\,y_2,\,y_3\big\}$.

If $d_{G_0}(z_1)\le 2$,
 then by a similar discussion as we made while comparing the $AZI(G_0)$ and $AZI(G_1)$,
we have  $AZI(G_0)>AZI(G_3)$, as
   $d_{G_0}(y_2)>d_{G_0}(y_1)$.

  Otherwise,  $d_{G_0}(w)\ge 3$ holds  for any  vertex $w\in N_{G_0}(y_1)\setminus V(C_g)$, then we also have $d_{G_0}(y_1)=3$ similarly with \eqref{32e}. With the similar reason as in Case 1, we can get another $k$-apex tree $G_4\in  \mathscr{U}_n^k$ such that $AZI(G_0)>AZI(G_4)$, a contradiction.    \end{proof}

 \begin{lemma}\label{37l} If  $k\ge 3$,  then the unique cycle of $G_0$ is a triangle.  \end{lemma}
\begin{proof}  By contradiction, we assume that $g\ge 4$. If $g\ge 5$, then by applying Lemma \ref{36l} step by step,  we can obtain a unicyclic graph $G^*\in  \mathscr{U}^k_n$ such that the unique cycle of $G^*$ is a $C_4$ and $AZI(G^*)=AZI(G_0)$. Thus, we may suppose that $G_0$ itself  is a unicyclic graph with unique cycle $C_4$   and $V(T_{w_0})$ is a $k$-quasi vertex set of  $G_0$  such that $V(C_4)=\{w_0,y_1,y_2,y_3\}$ with $N_{G_0}(y_1)\cap N_{G_0}(y_3)=\{w_0,y_2\}$. By Lemma \ref{31l}, each vertex of $C_4$ is a branching vertex.

 By the proof of Lemma \ref{36l} and the symmetry of $y_1$ and $y_3$, we have  $d_{G_0}(y_1)=d_{G_0}(y_2)=d_{G_0}(y_3)=s$. Let $G_1=G_0+y_1y_3+y_2z-y_1y_2-y_2y_3$, where $z$ is a pendant vertex of $T_{y_1}$ such that the distance between $z$ and $y_1$ is as small as possible. By Lemma \ref{31l}, $G_1\in \mathscr{U}_n^k$.

\par\noindent  {\bf Case 1. $s\ge 5$.} By the choice of $G_0$ and Proposition \ref{21p},  we can conclude that $s=5$ and   $T_{y_i}=K_{1,3}$ with   $y_i$ being   the center vertex of $K_{1,3}$ for $1\le i\le 3$. This implies that
\begin{align*}AZI(G_0)-AZI(G_1)=\Psi(5,5)+4\Psi(5,1)-16-3\Psi(4,1)>15,\end{align*}
a contradiction.

\par\noindent  {\bf Case 2. $s=4$.} By the choice of $G_0$ and Proposition \ref{21p},  we can conclude that either $T_{y_i}=K_{1,2}$ with   $y_i$ being   the center vertex of $K_{1,2}$ or $T_{y_i}=P_4$ with $y_i$ being a 2-vertex of $P_4$, where $i\in \{1,2,3\}$.  Since  $T_{y_2}\in \big\{K_{1,2},\,P_4\big\}$ and $zy_1\in E(G_0)$, we have \begin{align*}AZI(G_0)-AZI(G_1)\ge \Psi(4,4)+3\Psi(4,1)-16-2\Psi(3,1)>3 \,\,\text{by Lemma \ref{21l}}, \end{align*}
a contradiction.

\par\noindent  {\bf Case 3. $s=3$.} In this case, $d_{G_0}(y_1)=d_{G_0}(y_2)=d_{G_0}(y_3)=3$ and thus $|V(T_{y_2})|\ge 7$ by Lemma \ref{35l} and the choice of $G_0$. By Proposition \ref{21p}, we have  $T_{y_2}=Z_q$ and $y_2$ is the  special pendant vertex of $Z_q$, where  $q\ge 7$. Suppose that $w$ is the unique branching vertex of $V(Z_q)\setminus \{y_2\}$ and $z_1\in N_{G_0}(y_1)\setminus V(C_4)$.  Let   $G_2=G_0+w_0y_2+y_1w-w_0y_1-y_2y_1$. By Lemma \ref{31l}, $G_2$ is also a $k$-apex trees with a triangle as its  unique cycle.

 Since $|V(T_{y_1})|\ge k\ge 3$, we have  $d_{G_0}(z_1)\ge 2$ and thus
 \begin{align*}AZI(G_0)-AZI(G_2)&=\Psi(3,3)+\Psi\big(d_{G_0}(z_1),3\big)+(q-3)\Psi(q-2,1)-16-(q-3)\Psi(q-1,1)
 \\&\ge \Psi(3,3)-8>0,\end{align*}
by Lemma \ref{21l}, a contradiction.\end{proof}

We need to introduce the following two classes of graphs before proceeding further:\par\noindent
$(i)$  Let $D_{s_1,s_2,s_3}$ be the unicyclic graph with unique cycle $C_3=u_1u_2u_3$ such that $T_{u_1}=K_{1,s_1-1}$ with $u_1$ being the center vertex of $K_{1,s_1-1}$ and  $T_{u_i}=K_{1,s_i-1}$ with  $u_i$ being a  pendant vertex of $K_{1,s_i-1}$ for $i\in \{2,3\}$, where $s_1\ge 2$ and $\min\{s_2,s_3\}\ge 3$. \par\noindent
$(ii)$  Let $E_{s_1,s_2,s_3}$ be the unicyclic graph with unique cycle $C_3=u_1u_2u_3$ such that $T_{u_i}=K_{1,s_i-1}$ with $u_i$ being the center vertex of $K_{1,s_i-1}$ for $i\in \{1,2\}$  and  $T_{u_3}=Z_{s_3}$ with  $u_3$ being the special   pendant vertex of $Z_{s_3}$, where $\min\{s_1,s_2\}\ge 2$ and $s_3\ge 4$. \par
For instance, the graphs $D_{5,6,7}$ and $E_{5,6,7}$ are shown in Fig. \ref{21f}.
\par\medskip \noindent{\bf Proof of Proposition \ref{31p}: }   By Lemmas  \ref{31l} and  \ref{37l},     $C_3:\,=u_1u_2u_3$ is the unique cycle of $G_0$  with $$\text{$d_{G_0}(u_1)=\max\big\{d_{G_0}(u_i):\,1\le i\le 3\big\}$ and  $\min\big\{|V(T_{u_i})|:\,1\le i\le 3\big\}=k$.}$$ \par

\par\noindent
{\bf Case 1. $k\ge 7$.} By Proposition \ref{21p}, we have  $G_0=A_{k,s_1,s_2}$, where $s_1\ge s_2\ge k\ge 7$.
By Lemma \ref{32l}, we have  $G_0=A_{k,\left\lfloor\frac{n-k}{2}\right\rfloor, \left\lceil \frac{n-k}{2}\right\rceil}$, as desired.

\par\noindent
{\bf Case 2. $k=6$.}    Since   $\min\big\{|V(T_{u_i})|:\,1\le i\le 3\big\} =k=6$, by  Proposition \ref{21p} and  Lemma \ref{35l}, we have  $d_{G_0}(u_1)\ge 4$ and $d_{G_0}(u_2)=d_{G_0}(u_3)=3$. Since  $\min\big\{|V(T_{u_2})|, |V(T_{u_3})|\big\}\ge 6$, by  Proposition \ref{21p} and   Lemma   \ref{32l}, we have  $G=B_{6,\left\lceil\frac{n}{2}\right\rceil-3,\left\lfloor\frac{n}{2}\right\rfloor-3}$.

\par\noindent
{\bf Case 3. $k=5$.} If $d_{G_0}(u_1)=3$, then $d_{G_0}(u_2)=d_{G_0}(u_3)=3$. Since $\min\big\{|V(T_{u_i})|:\,1\le i\le 3\big\} =5$ and  Lemma \ref{35l},  there exists some $u_i$ such that   $T_{u_i}=K_{1,4}$ and  $u_i$ is the center  vertex of $K_{1,4}$, a contradiction. Thus,    $d_{G_0}(u_1)\ge 4$.

By $\min\big\{|V(T_{u_i})|:\,1\le i\le 3\big\} =k=5$, Lemma \ref{33l} and  Proposition \ref{21p}, we can conclude that  $d_{G_0}(u_2)=d_{G_0}(u_3)=3$  and either   $T_{u_i}=K_{1,4}$ with   $u_i$ being the pendant vertex of  $K_{1,4}$ or $T_{u_i}=Z_{q_i}$ with $u_i$ being the special pendant vertex of $Z_{q_i}$, where $q_i\ge 6$ and   $i\in \{2,3\}$. Once again,  $\min\big\{|V(T_{u_i})|:\,1\le i\le 3\big\} =5$  and  Proposition \ref{21p} imply that  $T_{u_1}=K_{1,d_{G_0}(u_1)-2}$ with $u_1$ as the center vertex of $K_{1,d_{G_0}(u_1)-2}$ and $6\le d_{G_0}(u_1)\le 7$. Next suppose that $|V(T_{u_2})|\ge |V(T_{u_3})|\ge 5$.

 If $|V(T_{u_3})|=5$, then $T_{u_3}=K_{1,4}$ with $u_3$ being a pendant vertex of $K_{1,4}$. Since $n\ge 18$,  $G_0=C_{d_{G_0}(u_1)-1,5,n-d_{G_0}(u_1)-4}$ by Proposition \ref{21p}, where  $n-d_{G_0}(u_1)-4\ge 7$ and $6\le d_{G_0}(u_1)\le 7$.

 If $|V(T_{u_3})|\ge 6$, then $G_0=B_{5,s_1,s_2}$ by Proposition \ref{21p}, as   $\min\big\{|V(T_{u_i})|:\,1\le i\le 3\big\} =k=5$. By  Lemma \ref{32l}, we have  $G_0=B_{5,\lceil\frac{n-5}{2}\rceil,\lfloor\frac{n-5}{2}\rfloor}$.

Since  $n\ge 18$, we have  $AZI(C_{5,5,n-10})<AZI(C_{6,5,n-11})$ by Lemma \ref{34l}. From the above arguments, we can conclude that   $G_0\in \left\{C_{5,5,n-10},\,B_{5,\lceil\frac{n-5}{2}\rceil,\lfloor\frac{n-5}{2}\rfloor}\right\}.$

By an elementary computation, we have
\begin{align*}AZI\left(B_{5,\lceil\frac{n-5}{2}\rceil,\lfloor\frac{n-5}{2}\rfloor}\right)-AZI(C_{5,5,n-10})
\ge&  16+\frac{(n-9)^3}{(n-11)^2} - 3\left(\frac{4}{3}\right)^3-\left(\frac{12}{5}\right)^3-\frac{(n-12)^3}{(n-13)^2}
\\>& \frac{(n-9)^3}{(n-11)^2} -\frac{(n-12)^3}{(n-13)^2}-5\\=&\frac{9n^3 - 332n^2 + 4053n - 16358}{(n-11)^2(n-13)^2}>0,\end{align*}
and thus   $G_0=C_{5,5,n-10}$.

\par\noindent
{\bf Case 4. $k=4$.} By the choice of $G_0$,  Lemma  \ref{35l} and $\min\big\{|V(T_{u_i})|:\,1\le i\le 3\big\} =k=4$, we have $d_{G_0}(u_1)\ge 4$. Suppose that $d_{G_0}(u_2)\ge d_{G_0}(u_3)$.

If $d_{G_0}(u_2)\ge 4$, then Proposition \ref{21p}   confirms  that  either   $T_{u_i}=K_{1,3}$ with $u_i$ being the center vertex of $K_{1,3}$ or   $T_{u_i}=P_4$ with $u_i$ being a 2-vertex of $P_4$ for $i\in \{1,2\}$. Combining this with $n\ge 15$, we have  $|V(T_{u_3})|\ge 7$. By Proposition \ref{21p},   $T_{u_3}=Z_q$ with $u_3$ being the special  pendant vertex of $Z_q$ for $q\ge 7$. We claim that $T_{u_i}=K_{1,3}$ with $u_i$ being the center vertex of $K_{1,3}$ for $i\in \{1,2\}$. Otherwise, assume that $T_{u_1}=P_4$ with $u_1$ being a 2-vertex of $P_4$. Let $G_1$ be the graph obtained from $G_0$ by replacing $T_{u_1}$ by $K_{1,3}$ such that $u_1$ is the center vertex of $K_{1,3}$. Suppose that $d_{G_0}(u_2)=s$. Then,   $4\le s\le 5$ and thus
\begin{align*}AZI(G_0)-AZI(G_1)&=\Psi(4,3)+\Psi(s,4)+\Psi(4,1)+16-\Psi(5,3)-\Psi(5,s)-3\Psi(5,1)\\&\ge \Psi(4,3)+\Psi(5,4)+\Psi(4,1)+16-\Psi(5,3)-\Psi(5,5)-3\Psi(5,1)>3.\end{align*}   This is impossible and thus $G_0=E_{4,4,n-8}$.

If $d_{G_0}(u_2)=3$, then $d_{G_0}(u_3)=3$. Since $d_{G_0}(u_1)\ge 4$,   Proposition \ref{21p} and  Lemma \ref{35l} imply that    $T_{u_1}=K_{1,d(u_1)-2}$ with $u_1$ being the center vertex of $K_{1,d(u_1)-2}$ for $5\le d_{G_0}(u_1)\le 7$.

When $6\le d_{G_0}(u_1)\le 7$, then we suppose that  $T_{u_2}=K_{1,3}$ with $u_2$ being a pendant vertex of $K_{1,3}$  by Lemma  \ref{33l}. By Lemmas \ref{33l}--\ref{34l}, Proposition \ref{21p}   and $n\ge 15$,   we have $G_0\in \big\{D_{6,4,5}, C_{5,4,n-9}\big\}$.

When $d_{G_0}(u_1)=5$, then Lemmas \ref{32l}--\ref{33l} and  Proposition \ref{21p}  imply that    $G_0\in \big\{C_{4,4,n-8}, C_{4,5,n-9}\big\}$ for  $n\ge 15$ or $G_0=B_{4,\lceil\frac{n-4}{2}\rceil, \lfloor\frac{n-4}{2}\rfloor}$ for  $n\ge 16$.

If $n\ge 15$, then
\begin{align*}&AZI\left(B_{4,\lceil\frac{n-4}{2}\rceil, \lfloor\frac{n-4}{2}\rfloor}\right)-AZI(C_{4,4,n-8})\ge 16+ \frac{(n-8)^3}{(n-10)^2}-\left(\frac{9}{4}\right)^3-2\left(\frac{3}{2}\right)^3
-\frac{(n-10)^3}{(n-11)^2}\\=&\frac{119n^4 - 4422n^3 + 60483n^2 - 359316n + 777372}{64(n-11)^2(n-10)^2}>0,\\ &AZI(C_{4,5,n-9})-AZI(C_{4,4,n-8})=\left(\frac{12}{5}\right)^3
+3\left(\frac{4}{3}\right)^3+\frac{(n-11)^3}{(n-12)^2}
-\left(\frac{9}{4}\right)^3-2\left(\frac{3}{2}\right)^3-\frac{(n-10)^3}{(n-11)^2}
\\[1mm]>& 2.75+\frac{(n-11)^3}{(n-12)^2}
-\frac{(n-10)^3}{(n-11)^2}
\\[1mm]=&\frac{7n^4 - 322n^3 + 5563n^2 - 42772n + 123460}{4(n-11)^2(n-12)^2}>0, \,\,\text{and}  \\[2mm] &AZI(E_{4,4,n-8})-AZI(C_{4,4,n-8})=\left(\frac{25}{8}\right)^3+3\left(\frac{5}{4}\right)^3-2\left(\frac{3}{2}\right)^3
-2\left(\frac{9}{4}\right)^3>6.5.\end{align*}
Combining this with $AZI(D_{6,4,5})>AZI(C_{4,4,7})$, we have    $G_0=C_{4,4,n-8}$. \qed

\section{$G_0$ is not a $c$-cyclic graph with $c\ge 4$}

In this section we discard, as AZI minimizers in $\mathbb{T}^k_n$, the $k$-apex trees whose cyclomatic number exceeds $3$.

\begin{proposition}\label{41p}
Let $G_0$ be a graph in $\mathbb{T}^k_n$ minimizing the AZI. Then $G_0$ is not a $c$-cyclic graph with $c\ge 4$.
\end{proposition}

\begin{lemma}\label{41l}If $n\ge 15$, then $$AZI\left(A_{5,5,n-10}\right)>AZI\left(B_{6,\left\lceil\frac{n}{2}\right\rceil-3,\left\lfloor\frac{n}{2}\right\rfloor-3}\right)>AZI(C_{5,5,n-10})>AZI(C_{4,4,n-8}).$$   \end{lemma}
\begin{proof} Since $\frac{6123}{64}-\frac{7088}{81}>8.15$, $\frac{7088}{81}-\frac{6961}{78}>-1.75$ and $\frac{6961}{78}-\frac{5289}{64}>6.5$, by  \eqref{42e}, \eqref{43e} and  \eqref{45e}, we have
\begin{eqnarray*}
% \nonumber to remove numbering (before each equation)
  AZI(A_{5,5,n-10})-AZI\left(B_{6,\left\lceil\frac{n}{2}\right\rceil-3,\left\lfloor\frac{n}{2}\right\rfloor-3}\right) &>& \frac{f_1(n)}{20(n-11)^2(n-13)^2}>0, \\
  AZI\left(B_{6,\left\lceil\frac{n}{2}\right\rceil-3,\left\lfloor\frac{n}{2}\right\rfloor-3}\right)-AZI(C_{5,5,n-10}) &>&  \frac{f_2(n)}{4(n-12)^2(n-13)^2}>0,\\
 \text{and}\,\, AZI(C_{5,5,n-10}) -AZI(C_{4,4,n-8}) &>& \frac{f_3(n)}{2(n-11)^2(n-13)^2}>0,
\end{eqnarray*}
where $f_1(n)=83n^4 - 4164n^3 + 77946n^2 - 645332n + 1993947,$  $f_2(n)=9n^4 - 414n^3 + 7081n^2 - 53320n + 148976$ and $f_3(n)= 9n^4 - 432n^3 + 7770n^2 - 62056n + 185661.$  \end{proof}

\begin{lemma}\label{40l} Let $G_0$ be a graph in $\mathbb{T}^k_n$ minimizing the AZI with $k\ge 4$.  If $n\ge 15$, then \begin{align}\label{liun2}AZI\left(G_0\right)<2\left(\frac{12}{5}\right)^3+64+(n-9)\left(\frac{n-1}{n-2}\right)^3.\end{align} Furthermore, if $AZI(G_0)=AZI(C_{4,4,n-8})$, then \begin{align}\label{liun3}AZI(G_0)<3\left(\frac{9}{4}\right)^3+48+(n-8)\left(\frac{n-1}{n-2}\right)^3.\end{align}   \end{lemma}
\begin{proof}  Assume that \eqref{liun3} does not hold. Then,
\begin{align}\label{53e}AZI(G_0)-AZI(C_{4,4,n-8})\ge& 3\left(\frac{9}{4}\right)^3+48+(n-8)\left(\frac{n-1}{n-2}\right)^3-\frac{5289}{64}-\frac{(n-10)^3}{(n-11)^2}
 \\\nonumber=&\frac{3f_1(n)}{32(n-11)^2 (n-2)^3}>0\,\,\text{for $n\ge 15$},\end{align}
 where $f_1(n)=27n^5 - 948n^4 + 11667n^3 - 59170n^2 + 111316n - 70168.$  This implies that \eqref{liun3} holds. Next, assume that \eqref{liun2} does not hold, that is,
\begin{align} \label{51e}AZI(G_0)\ge  2\left(\frac{12}{5}\right)^3+64+(n-9)\left(\frac{n-1}{n-2}\right)^3.\end{align}
Since $2\left(\frac{12}{5}\right)^3+64-\frac{6259}{64}>-6.15$ and $2\left(\frac{12}{5}\right)^3+64-\frac{7088}{81}>4$, we have
 \begin{align*}\nonumber AZI(G_0)-AZI(A_{7,7,n-14})> &\frac{f_2(n)}{20(n-2)^3(n-17)^2}>0\,\,\text{ for $n\ge 21$},\,\,\text{and}\,\,\\[2mm]AZI(G_0)-AZI\left(B_{6,\lceil\frac{n}{2}\rceil-3,\lfloor\frac{n}{2}\rfloor-3}\right)> &\frac{f_3(n)}{(n-2)^3(n-11)^2(n-13)^2}>0\,\,\text{ for $n\ge 18$},  \end{align*}   where
 $f_2(n)=37n^5 - 1900n^4 + 33325n^3 - 226450n^2 + 477940n - 318964$ and $f_3(n)=4n^7 - 246n^6 + 6060n^5 - 76199n^4 + 514188n^3 - 1786204n^2 + 2805268n - 1607335$.
 By \eqref{51e}, we also have    $AZI(G_0)>AZI(C_{5,5,n-10})$ for $15\le n\le 17$. This is contrary with the choice of $G_0$,  Corollary \ref{31c} and Lemma \ref{41l}.
\end{proof}

 \begin{lemma}\label{42l}{\em\cite{MH1}}
Let $G$ be a $c$-cyclic  graph with $n$ vertices.  If $n\geq 2c+1$ and $c\geq 3$, then $AZI(G)\geq AZI(W_{n;c,1})$, where the  equality holds  if and only if  $G=W_{n;c,1}$.
 \end{lemma}
 \begin{lemma}\label{43l}{\em\cite{Xiaod}}
Let $G$ be a connected graph with two non-adjacent vertices $u$ and $v$. If $\alpha \leq 1/2$ and $\alpha \neq 0$, then $ABC_\alpha(G + uv) > ABC_\alpha(G)$.
 \end{lemma}

\begin{lemma}\label{44l} Let $n\ge 15$, then $AZI(W_{n;5,1})>AZI(A_{5,5,n-10})$. \end{lemma}
\begin{proof}
 By   \eqref{41e} and \eqref{42e},
 \begin{align*}AZI(W_{n;5,1})-AZI(A_{5,5,n-10}) & =  (n-8)\left(\frac{n-2}{n-3}\right)^3+96-\frac{6123}{64}-\frac{(n-12)^3}{(n-13)^2}\\
 &> (n-8)\left(\frac{n-2}{n-3}\right)^3-\frac{(n-12)^3}{(n-13)^2}\\
 &= \frac{5f(n)}{(n - 3)^3(n-13)^2}>0,
 \end{align*}
 where $f(n)=n^5 - 38n^4 + 517n^3 - 3056n^2 + 7816n - 7168$.   \end{proof}

  In the following, let $G_0$ be the extremal graph with the minimum AZI among all $k$-apex trees with $n\ge 3k+3$ vertices, where $k\ge 4$.

If $\delta\ge 2$, then $c\ge 3$, otherwise, $G_0$ contains at most two 2-vertices in different cycles of $G_0$ such that the deletion of them   will be resulted in a tree, contrary with $k\ge 4$. Since $n\ge 15$, by Lemma \ref{43l}, Corollary \ref{21c}   and \eqref{42e}, \begin{align*} AZI(G_0)-AZI(A_{5,5,n-10})&\ge 8(n+2)-\frac{6123}{64}-\frac{(n-12)^3}{(n-13)^2}\\ \nonumber &>8n-80-\frac{(n-12)^3}{(n-13)^2}=\frac{7n^3 - 252n^2 + 3000n - 11792}{(n-3)^2}>0,\end{align*}
a contradiction. Combining the above fact with Proposition  \ref{31p},  we may assume that $$\text{   $c\ge 2$ and   $\delta=1$.}$$

As in \cite{MH21}, a {\bf special cycle} is a cycle with each vertex being a branching vertex. Recall  that $G_0$ is a $c$-cyclic graph with  exactly $p$ pendant vertices. The following result implies that  $p\ge n-11+c$.

\begin{lemma}\label{45l}   If $2\le c\le 4$, then   $p\ge n-11+c$. Moreover,  if either $k\ge 5$ or $2\le c\le 3$, then $G_0$ contains a special cycle.    \end{lemma}
\begin{proof} Assume that $p\le n-12+c$.
\par\noindent {\bf Case 1. $k=4$.} By Corollary \ref{21c}  and  \eqref{43e}, we have
\begin{align}\nonumber AZI(G_0)-AZI(C_{4,4,n-8})\ge & \, 8(n+c-1-p)+p\left(\frac{n-1}{n-2}\right)^3-\frac{5289}{64}-\frac{(n-10)^3}{(n-11)^2}\\
\nonumber \ge& \; 88+(n-12+c)\left(\frac{n-1}{n-2}\right)^3-\frac{5289}{64}-\frac{(n-10)^3}{(n-11)^2}\\
 \label{46e}\ge & \; 88+(n-10)\left(\frac{n-1}{n-2}\right)^3-\frac{5289}{64}-\frac{(n-10)^3}{(n-11)^2}\end{align}\begin{align}
\nonumber>& \; 5+(n-10)\left(\frac{n-1}{n-2}\right)^3-\frac{(n-10)^3}{(n-11)^2}\\
\nonumber =& \; \frac{f_1(n)}{(n-2)^3(n-11)^2}>0\,\,\text{ for $n\ge 15$},
\end{align}
where $f_1(n)=6n^5 - 192n^4 + 2163n^3 - 10145n^2 + 18569n - 11630$.

\par\noindent {\bf Case 2. $k\ge 5$.} We first show that $G_0$ contains at least one special cycle. Otherwise, assume that $G_0$ contains no special cycle  and thus each cycle of $G_0$ contains a 2-vertex. Let $G_1=G_0-u$, where $u$ is a 2-vertex of some cycle of $G_0$. Then, $G_1$ is a $(c-1)$-cyclic graph with $n-1$ vertices. By our hypothesis, every  cycle of $G_1$ is not a special cycle and thus defined $G_2=G_1-v$, where $v$ is a 2-vertex of some cycle of $G_2$. Now, $G_2$ is a  $(c-2)$-cyclic graph    with $n-2$ vertices and $G_2$ contains no special cycle. With the similar reason, since $2\le c\le 4$, we can obtain a tree by deleting at most $c$ vertices from $G_0$, which implies that $G_0$ is not a $k$-apex tree with $k\ge 5$. Thus, $G_0$ contains at least one special cycle, which implies that   $G$ contains at least three $\left(3^{\ge},3^{\ge }\right)$ edges. Combining this with Corollary \ref{21c}, $p\le n-12+c$ and  \eqref{42e}, we have   \begin{align*}AZI(G_0)-AZI(A_{5,5,n-10})\ge & 3\left(\frac{9}{4}\right)^3+64+(n-10)\left(\frac{n-1}{n-2}\right)^3-\frac{6123}{64}-\frac{(n-12)^3}{(n-13)^2}
\\=&\frac{f_2(n)}{2(n-2)^3(n-13)^2}>0\,\,\text{ for $n\ge 15$}, \end{align*}   where
  $f_2(n)=11n^5 - 400n^4 + 5033n^3 - 25564n^2 + 48566n - 31028$.    \end{proof}

Hereafter, let $G^*_0$ be the graph obtained from $G_0$ by deleting all the pendant vertices of $G_0$. Clearly, if $G_0$ is a $c$-cyclic graph with $p$ pendant vertices and $n$ vertices, then $G^*_0$ is a $c$-cyclic graph with $n-p$ vertices.
\par\medskip \noindent{\bf Proof of Proposition \ref{41p}:} Let $G_0$ be a $c$-cyclic graph. In view of Lemmas \ref{42l}, \ref{43l} and \ref{44l}, for $n\geq 15$ and $c\ge 5$, we have the following chain of inequalities
\[AZI(G)\geq AZI(W_{n;c,1}) \geq  AZI(W_{n;5,1})>AZI(A_{5,5,n-10}) .\]
Therefore, if $G_0\in\mathbb{T}^k_n$ minimizes the AZI and $n\geq15$, then $c\le 4$.

Next, we shall prove that it is impossible that  $G_0$ is a $4$-cyclic graph. By Lemma \ref{45l},   we have $n-7\le p\le n-5$.

\par\noindent{\bf Case 1. $p=n-7.$} If $G_0$ contains three $\left(3^{\ge},3^{\ge }\right)$ edges,  then  \begin{align*}AZI(G_0)\ge 3\left(\frac{9}{4}\right)^3+56+(n-7)\left(\frac{n-1}{n-2}\right)^3.\end{align*}
 Since $3\left(\frac{9}{4}\right)^3+58>2\left(\frac{12}{5}\right)^3+64$, the result follows from \eqref{51e}. Thus,  we suppose that $G_0$ contains at most two $\left(3^{\ge},3^{\ge }\right)$ edges and thus $k=4$ by Lemma \ref{45l}.

 If $G_0$ contains one $\left(3^{\ge},3^{\ge }\right)$ edge, then the result follows from \eqref{53e}, as $\left(\frac{9}{4}\right)^3+72>3\left(\frac{9}{4}\right)^3+48$. Thus, $G_0$ contains no $\left(3^{\ge},3^{\ge }\right)$ edge. Combining this with $p=n-7$ and $c=4$, we can conclude that      $G_0$ must be obtained from $K_{2,5}$ by attaching $r$ and $n-7-r$ pendant vertices to   the two 5-vertices of $K_{2,5}$, respectively. Combining this with $k=4$, we have $r\ge 3$ and thus  $\Delta\le n-5$.   By Corollary \ref{21c} and \eqref{43e}, we have \begin{align*}AZI(G_0)-AZI(C_{4,4,n-8})&\ge 80+(n-7)\left(\frac{n-5}{n-6}\right)^3-\frac{5289}{64}-\frac{(n-10)^3}{(n-11)^2}
 \\&>(n-7)\left(\frac{n-5}{n-6}\right)^3-\frac{(n-10)^3}{(n-11)^2}-3\\&=\frac{f(n)}{(n-11)^2(n-6)^3}>0 \,\,\text{for $n\ge 15$},\end{align*}  where $f(n)=n^5 - 43n^4 + 709n^3 - 5615n^2 + 21440n - 31717$.

  \par\noindent{\bf Case 2. $p=n-6.$} Since $p=n-6$, $G^*_0$ is also a 4-cyclic graphs with   at most six vertices. From the table of  connected graphs with at most six vertices (see   \cite{Dc1}), we can see that  $G^*_0$ contains six   $\left(3^{\ge},3^{\ge }\right)$ edges, or four  $\left(4^{\ge},3^{\ge }\right)$ edges, or three $\left(4^{\ge},4^{\ge }\right)$ edges,  or one  $\left(4^{\ge},3^{\ge }\right)$ edge and two $\left(5^{\ge},3^{\ge }\right)$ edges, or one $\left(5^{\ge},5^{\ge }\right)$ edge.   By Corollary \ref{21c},  \begin{align*}AZI(G_0)\ge 6\left(\frac{9}{4}\right)^3+24+(n-6)\left(\frac{n-1}{n-2}\right)^3.\end{align*} Since $6\left(\frac{9}{4}\right)^3+24>2\left(\frac{12}{5}\right)^3+64$, the result follows from \eqref{51e}.

 \par\noindent{\bf Case 3. $p=n-5.$} Since $p=n-5$, $G^*_0$ is also a 4-cyclic graphs with   at most five vertices. From the table of  connected graphs with   five  vertices (see   \cite{Dc1}), we can see that  $G^*_0$ contains one  $\left(3^{\ge},3^{\ge }\right)$ edge and four  $\left(4^{\ge},3^{\ge }\right)$ edges. By Corollary \ref{21c},
 \begin{align*}AZI(G_0)\ge \left(\frac{9}{4}\right)^3+4\left(\frac{12}{5}\right)^3+24+(n-5)\left(\frac{n-1}{n-2}\right)^3.\end{align*}
Since $\left(\frac{9}{4}\right)^3+4\left(\frac{12}{5}\right)^3+24+4>2\left(\frac{12}{5}\right)^3+64$, the result follows from \eqref{51e}.   \qed

\section{$G_0$ is   not a bicyclic graph}
For $q\ge2$, let $Q_{s_1,s_2,\ldots,s_q}$ be a $(q-1)$-cyclic graph obtained from $q$ non-trivial paths $P_{s_1}, P_{s_2}, \ldots, P_{s_q}$ by firstly identifying the vertices  $u_1,u_2,\ldots,u_q$ and then identifying the vertices $v_1,v_2, \ldots, v_q$, where  $u_i,v_i$ are the two end vertices of $P_{s_i}$ for $i=1, 2,\ldots, q$, and $s_i=2$ for at most one $i\in\{1,2,\ldots,q\}$.
 In some papers,  $Q_{s_1,s_2,s_3}$ is also referred as a $\Theta$-graph.   From the definition, $Q_{s_1,s_2,\ldots,s_q}$ contains exactly two $q$-vertices and each of  the other vertices is a 2-vertex.   Since $G_0$ is a $k$-apex tree with $k\ge 4$, if $\mathcal {R}(G_0)=Q_{s_1,s_2,s_3}$, then the two 3-vertices of $Q_{s_1,s_2,s_3}$ are two vertices with degree at least four  in $G_0$.
Here is the main result of this section.
\begin{proposition}\label{51p}
Let $G_0$ be a graph in $\mathbb{T}^k_n$ minimizing the AZI. Then  $G_0$ is not   bicyclic.
\end{proposition}

We assume that $G_0$ is a bicyclic graph in the rest of this section.  By Lemma \ref{45l}, we can confirm that $p\ge n-9$ and $G_0$ contains a special cycle $C_s$.  Since $G_0$ is a bicyclic graph,   $G_0$ contains another cycle $C_t$.

\begin{lemma}\label{51l}  If $G_0$ is a bicyclic graph, then   $p\ge n-8$.  \end{lemma}
\begin{proof} By Lemma \ref{45l}, assume that $p=n-9$. If    $G_0$ contains  two  $\big(4^{\ge},3^{\ge}\big)$ edges or four  $\big(3^{\ge},3^{\ge}\big)$ edges, then    Corollary \ref{21c}  implies that \eqref{51e} holds, contrary with Lemma \ref{40l}. Thus,    $G_0$ contains at most three $\big(3^{\ge},3^{\ge}\big)$ edges and $G_0$ contains  at most one $\big(4^{\ge},3^{\ge}\big)$ edge. This implies that  the special cycle of $G_0$ is a triangle, each vertex of which is a 3-vertex, and thus $\mathcal {R}(G_0)=B_1$, where $x_0y_0\not\in E(G_0)$ (see Fig. \ref{21f}).

 Since $G_0$ contains at most three  $\big(3^{\ge},3^{\ge}\big)$ edges, the two neighbor of $y_0$ in $C_t$ are two 2-vertices of $G_0$. Combining this with $p=n-9$ and $G_0$ being a $k$-apex tree with $k\ge 4$, we can conclude that $C_t$ is also a triangle and  the two 3-vertices  of $C_s$  are adjacent with one 2-vertex of $G_0$, which is adjacent to a pendant vertex, respectively. This implies that  $k=4$ and $G_0$ has three  $\big(3^{\ge},3^{\ge}\big)$ edges and night  edges with a 2-vertex as their end vertex.
 Since $8\times 8+3\left(\frac{9}{4}\right)^3>88$, we have $AZI(G_0)>AZI(C_{4,4,n-8})$ by \eqref{46e}.
 \end{proof}

\begin{lemma}\label{52l} Let $G_0$ be a bicyclic graph. If one of the followings holds, then $p\ge   n-6$:\par\noindent
   $(1)$ $G_0$ contains    eight $\left(3^{\ge},3^{\ge }\right)$ edges, \par\noindent
    $(2)$ $G_0$ contains    one   $\left(4^{\ge},4^{\ge }\right)$ edge, two  $\left(4^{\ge},3^{\ge }\right)$ edges and one $\left(3^{\ge},3^{\ge }\right)$,  \par\noindent
  $(3)$ $G_0$ contains  four  $\left(4^{\ge},3^{\ge }\right)$ edges and  one  $\left(3^{\ge},3^{\ge }\right)$ edges.   \end{lemma}
\begin{proof}  Since $8\left(\frac{9}{4}\right)^3+2>4\left(\frac{12}{5}\right)^3+\left(\frac{9}{4}\right)^3+24+2>2\left(\frac{12}{5}\right)^3+\left(\frac{16}{6}\right)^3+\left(\frac{9}{4}\right)^3+32+2
>2\left(\frac{12}{5}\right)^3+64$,    result follows from \eqref{51e}.  \end{proof}

Similarly with Lemma  \ref{52l}, we can show the following two lemmas:
\begin{lemma}\label{53l}   Let $G_0$ be a bicyclic graph. If one of the followings holds, then $p\ge   n-5$:\par\noindent
  $(1)$ $G_0$ contains  one   $\left(4^{\ge},4^{\ge }\right)$ edge, two  $\left(4^{\ge},3^{\ge }\right)$ edges and three   $\left(3^{\ge},3^{\ge }\right)$ edges, \par\noindent
  $(2)$ $G_0$ contains  one $\left(4^{\ge},4^{\ge }\right)$ edge  and four    $\left(4^{\ge},3^{\ge }\right)$ edges, \par\noindent
  $(3)$ $G_0$ contains  three $\left(4^{\ge},4^{\ge }\right)$ edges,\par\noindent
  $(4)$ $G_0$ contains   four $\left(5^{\ge},3^{\ge }\right)$ edges and one   $\left(3^{\ge},3^{\ge }\right)$ edge,\par\noindent
  $(5)$ $G_0$ contains  four  $\left(4^{\ge},3^{\ge }\right)$ edges  and three   $\left(3^{\ge},3^{\ge }\right)$ edges.
  \end{lemma}

\begin{lemma}\label{54l}   Let $G_0$ be a bicyclic graph. If one of the followings holds, then $p=n-4$:\par\noindent   $(1)$ $G_0$ contains  one  $\left(5^{\ge},4^{\ge }\right)$ edges, four  $\left(4^{\ge},3^{\ge }\right)$ edges and one   $\left(3^{\ge},3^{\ge }\right)$ edge,\par\noindent
  $(2)$ $G_0$ contains two $\left(5^{\ge},4^{\ge }\right)$ edges, one   $\left(4^{\ge},3^{\ge }\right)$ edge and one $\left(3^{\ge},3^{\ge }\right)$ edge, \par\noindent
  $(3)$ $G_0$ contains  one    $\left(6^{\ge},4^{\ge }\right)$ edge, one  $\left(6^{\ge},3^{\ge }\right)$ edge and two $\left(4^{\ge},3^{\ge }\right)$ edges.  \end{lemma}

   \begin{lemma}\label{55l}{\em \cite{MH21}}  If $G\in \mathbb{T}^3_n$ and $n\ge 14$, then $AZI(G)\ge \frac{8432}{125}+(n-8)\left(\frac{n-4}{n-5}\right)^3$.  \end{lemma}

   \begin{lemma}\label{56l}    If $k\ge 5$,  then no cycle of $G_0$ contains $2$-vertex.  \end{lemma}
\begin{proof} By contradiction, assume that $u_0v_0$ is an edge of some cycle of $G_0$ such that $u_0$ is a  2-vertex of $G_0$. By Lemma \ref{23l}, we have $G_0-u_0v_0\in \mathscr{U}^r_n$   with $k-1\le r\le k$. Combining this with   Lemma \ref{24l} and the choice of $G_0$, we may suppose that  \begin{align}\label{52e}\text{$G_0-u_0v_0\in \mathscr{U}^{k-1}_n$ and   $AZI(G_0)\ge AZI(G_0-u_0v_0)+8$.}\end{align}

\par\noindent{\bf Case 1. $k\ge 8$.} By   \eqref{52e}, Corollary \ref{31c} and Proposition \ref{31p}, we have  $$AZI(G_0)>AZI\left(A_{k-1,\left\lfloor\frac{n-k+1}{2}\right\rfloor, \left\lceil \frac{n-k+1}{2}\right\rceil}\right)>AZI\left(A_{k,\left\lfloor\frac{n-k}{2}\right\rfloor, \left\lceil \frac{n-k}{2}\right\rceil}\right),$$ and thus  the result already holds.

\par\noindent{\bf Case 2. $6\le k\le 7$.} By   \eqref{52e}, Lemma  \ref{41l} and Proposition \ref{31p}, we have  $$AZI(G_0)-
AZI\left(A_{5,5,n-10}\right)\ge AZI\left(C_{5,5,n-10}\right)+8-AZI\left(A_{5,5,n-10}\right)>1.5,$$ and thus  the result already holds.

\par\noindent{\bf Case 3. $k=5$.} By \eqref{52e},  $G_0-u_0v_0\in \mathscr{U}^4_n$. Combining this with  Proposition \ref{31p}, we have   \begin{align*}AZI(G_0)-AZI(C_{5,5,n-10})&\ge AZI\left(C_{4,4,n-8}\right)+8-AZI(C_{5,5,n-10})\\&>\frac{(n-10)^3}{(n-11)^2}-\frac{(n-12)^3}{(n-13)^2}
\\&=\frac{2(n^4 - 48n^3 + 859n^2 - 6794n + 20044)}{(n-11)^2\,(n-13)^2}>0\,\,\text{for $n\ge 15$}.\end{align*} This completes the proof of this result.   \end{proof}

 \begin{lemma}\label{57l}  Any cycle of    $ G_0$ contains no  $2$-vertex of $G_0$.   \end{lemma}
\begin{proof}  By contradiction, assume that $u_0v_0$ is an edge of some cycle of $G_0$ such that $u_0$ is a  2-vertex of $G_0$. By Lemma \ref{56l}, we have $k=4$ and thus   $G_0-u_0v_0\in \mathscr{U}^r_n$   with $3\le r\le 4$ by Lemma \ref{23l}. Combining this with   Lemma \ref{24l} and the choice of $G_0$, we may suppose that  $G_0-u_0v_0\in \mathscr{U}^3_n$.

 Since $C_s$ is a special cycle of $G_0$,    $G_0$  contains three $\left(3^{\ge},3^{\ge }\right)$ edges.
 If $p\le n-8$, then  Corollary \ref{21c} implies that  \eqref{53e} holds, contrary Lemma \ref{40l}.   Thus,  $n-4\le p\le n-7$.

  If $\mathcal {R}(G_0)=Q_{s_1,s_2,s_3}$, then let $x_0,y_0$ be the two 3-vertices of $Q_{s_1,s_2,s_3}$. Since  $G_0$ is a 4-apex tree, we can conclude that $|V(T_w)|\ge 3$ for each vertex $w\in \{x_0,y_0\}$, where  either $d_{G_0}(w)\ge 5$ or $w$ is adjacent with one non-pendant vertex. For convenience, we suppose that $s=s_1+s_2-2$.

 If $\mathcal {R}(G_0)\in \big\{B_1,\,B_2\big\}$, then let $x_0$ and $y_0$ be the unique cut vertex of $C_s$ and $C_t$, respectively, where $x_0=y_0$ when $\mathcal {R}(G_0)=B_2$ (see Fig. \ref{21f}).   Since $C_t$ contains a 2-vertex and $G_0$ is a 4-apex tree,   $|V(T_w)|\ge 3$ for each vertex $V(C_s)\setminus \{x_0\}$. This implies that either $d_{G_0}(w)\ge 4$ or $w$ is adjacent with one non-pendant vertex.

 \par\noindent{\bf Case  1.  $p=n-7$.}  If $G_0$ contains  three  $\left(4^{\ge},3^{\ge }\right)$ edges, then  the result follows from \eqref{53e}, as $3\left(\frac{12}{5}\right)^3+40+1>3\left(\frac{9}{4}\right)^3+48$. Thus, we may suppose that $G_0$ contains at most two  $\left(4^{\ge},3^{\ge }\right)$ edges. This implies that   $\mathcal {R}(G_0)\ne Q_{s_1,s_2,s_3}$.

 \par\noindent{\bf Subcase  1.1. $\mathcal {R}(G_0)=B_2$.} If $s\ge 4$, then $G_0$ contains    two $\left(4^{\ge},3^{\ge }\right)$ edges and two  $\left(3^{\ge},3^{\ge }\right)$ edges and thus the result follows from \eqref{53e}, as  $2\left(\frac{12}{5}\right)^3+2\left(\frac{9}{4}\right)^3+32>3\left(\frac{9}{4}\right)^3+48$. Thus, we have $s=3$ and   and the two neighbors of $C_t$ of $x_0$ are two 2-vertices of $G_0$.

  If $t=3$, then   $C_t$ contains a $(2,2)$-edge, say $u_1u_2$,  and thus \begin{align}\label{54e} AZI(G_0)\ge AZI(G_0-u_1u_2)+24-2\left(\frac{4}{3}\right)^3>AZI(G_0-u_1u_2)+19\end{align} by \eqref{21e}, as $d_{G_0}(x_0)\ge 4$. Since $19+\frac{8432}{125}>3\left(\frac{9}{4}\right)^3+48$, the result follows from \eqref{53e} and Lemma \ref{55l}. Combining this with $G_0$ being  a 4-apex tree, we have  $t\ge 4$. Since $p=n-7$, $C_t$ contains a 2-vertex of $G_0$ and $G_0$ is a 4-apex tree, we can conclude that   at least one    vertex of $V(C_s)\setminus \{x_0\}$ has  degree at least four and thus $G_0$ contains  three $\left(4^{\ge},3^{\ge }\right)$ edges, a contradiction.

 \par\noindent{\bf Subcase  1.2.  $\mathcal {R}(G_0)=B_1$.} If $\max\{s,t\}=4$, then each vertex of $V(C_s)\setminus \{x_0\}$ is a vertex with degree at least four, as $C_t$ contains a 2-vertex of $G_0$ and $G_0$ is a 4-apex tree. This implies that  $G_0$ contains three $\left(4^{\ge},3^{\ge }\right)$ edges, a contradiction. Otherwise, $s=t=3$.

 Since $G_0$ is a 4-apex tree and $G_0$ contains at most two $\left(4^{\ge},3^{\ge }\right)$ edges,   we can conclude that   $V(C_s)\setminus \{x_0\}$ has one vertex of degree at least four and exactly one 3-vertex, which is adjacent to one non-pendant vertex. Combining this with $p=n-7$, we have $x_0y_0\in E(G_0)$. Thus, $G_0$ contains  two $\left(4^{\ge},3^{\ge }\right)$ edges and two  $\left(3^{\ge},3^{\ge }\right)$ edges. Now,  the result follows from \eqref{53e}.

 \par\noindent{\bf Case 2. $p=n-6$.}

 \par\noindent{\bf Subcase 2.1. $\mathcal {R}(G_0)= Q_{s_1,s_2,s_3}$.}

  We first suppose that $x_0y_0\not\in E(G_0)$. Since $G_0$ is a $4$-apex tree and $p=n-6$, we have $\max\big\{d_{G_0}(x_0), d_{G_0}(y_0)\big\}\ge 6$ and $\min\big\{d_{G_0}(x_0), d_{G_0}(y_0)\big\}\ge 4$, and thus $G_0$ contains two $\left(6^{\ge},3^{\ge }\right)$ edges and two  $\left(4^{\ge},3^{\ge }\right)$ edges. Since $2\left(\frac{18}{7}\right)^3+ 2\left(\frac{12}{5}\right)^3+24>3\left(\frac{9}{4}\right)^3+48$, the result follows from \eqref{53e}.

 We second suppose that $x_0y_0\in E(G_0)$.  If $\max\big\{d_{G_0}(x_0), d_{G_0}(y_0)\big\}\ge 5$, then $G_0$ contains one $\left(5^{\ge},4^{\ge }\right)$ edge and two  $\left(4^{\ge},3^{\ge }\right)$ edges and thus the result follows from \eqref{53e}, as  $2\left(\frac{12}{5}\right)^3+ \left(\frac{20}{7}\right)^3+32>3\left(\frac{9}{4}\right)^3+48$. Thus, $x_0$ and $y_0$ are two $4$-vertices of $G_0$, each is adjacent with one branching  vertex, respectively, as $p=n-6$ and $G_0$ is a 4-apex tree. Now, \begin{align}\label{55e}\text{ $G_0$ contains one  $\left(4^{\ge},4^{\ge }\right)$ edge, two   $\left(4^{\ge},3^{\ge }\right)$ edges and one  $\left(3^{\ge},3^{\ge }\right)$ edge}.\end{align}  Since $2\left(\frac{12}{5}\right)^3+\left(\frac{16}{6}\right)^3+\left(\frac{9}{4}\right)^3+ 24+2>3\left(\frac{9}{4}\right)^3+48$,  the result follows from \eqref{53e}.

 \par\noindent{\bf Subcase 2.2. $\mathcal {R}(G_0)=B_1$.} Since $G_0$ is a 4-apex tree and $C_t$ contains a 2-vertex of $G_0$,  each vertex of $V(C_s)\setminus \{x_0\}$ has degree at least four. Thus,   $G_0$ contains one  $\left(4^{\ge},4^{\ge }\right)$ edge, two   $\left(4^{\ge},3^{\ge }\right)$ edges and one  $\left(3^{\ge},3^{\ge }\right)$ edge, and so the result follows from \eqref{55e}.

  \par\noindent{\bf Subcase 2.3. $\mathcal {R}(G_0)=B_2$.}   If  $\max\{s,t\}=4$, then each vertex of $C_s$ has degree at least four, as $G_0$ is a $4$-apex tree and $C_t$ contains a 2-vertex of $G_0$.  This implies that $G_0$ contains three $\left(4^{\ge},4^{\ge }\right)$ edges and thus the result follows from Lemma \ref{53l}$(3)$.

 Otherwise,   $s=t=3$. Since $G_0$ is a $4$-apex tree and $C_t$ contains a 2-vertex of $G_0$, $C_s$ contains at least two vertices with  degree at least four. By \eqref{55e}, we may suppose that the two vertices of $V(C_t)\setminus \{x_0\}$ are both 2-vertices of $G_0$. From \eqref{54e}, the result also holds.

 \par\noindent{\bf Case 3. $p=n-5$.} If  $\mathcal {R}(G_0)=B_2$, then   $s=t=3$ and the three vertices of $C_s$ has degree at least four, as $G_0$ is a $4$-apex tree. By Lemma \ref{53l}$(3)$, the result also holds.

 If $\mathcal {R}(G_0)=Q_{s_1,s_2,s_3}$, then   $\max\big\{d_{G_0}(x_0),\, d_{G_0}(y_0)\big\}\ge 6$, as $p=n-5$ and $G_0$ is a 4-apex tree. Thus,  $G_0$ contains one $\left(6^{\ge},4^{\ge }\right)$ edges, one   $\left(4^{\ge},3^{\ge }\right)$ edge and one   $\left(6^{\ge},3^{\ge }\right)$ edge. Since $\left(\frac{12}{5}\right)^3+ \left(\frac{18}{7}\right)^3+ \left(\frac{24}{8}\right)^3+27>3\left(\frac{9}{4}\right)^3+48$,   the result follows from \eqref{53e}.

 \par\noindent{\bf Case 4. $p=n-4$.} Then,  $\mathcal {R}(G_0)=Q_{s_1,s_2,s_3}$. Since  $G_0$ is a 4-apex tree and $p=n-4$,   $\min\big\{d_{G_0}(x_0),\, d_{G_0}(y_0)\big\}\ge 6$, and thus $G_0$ contains one $\left(6^{\ge},6^{\ge }\right)$ edges    and two   $\left(6^{\ge},3^{\ge }\right)$ edges. Since $2\left(\frac{18}{7}\right)^3+ \left(\frac{36}{10}\right)^3+16>3\left(\frac{9}{4}\right)^3+48$,   the result follows from \eqref{53e}.
  \end{proof}

 \begin{lemma}\label{58l}   $\mathcal {R}(G_0)$ contains no  cut vertex.   \end{lemma}
\begin{proof} We assume that $\mathcal {R}(G_0)$ contains a  cut vertex. By Lemma \ref{57l}, both $C_s$ and  $C_t$ are special cycles  and each of them contains a cut vertex of  $G_0$.
 By Lemma \ref{51l}, $p\ge n-8$.   Suppose that $V(C_s)\setminus \{x_0\}=\{w_1,w_2,\ldots,w_{s-1}\}$ and $V(C_t)\setminus \{y_0\}=\{z_1,z_2,\ldots,z_{t-1}\}$, where $x_0$ and $y_0$ are the cut vertex of $C_s$ and $C_t$, respectively and $x_0=y_0$ when $\mathcal {R}(G_0)=B_2$. For convenience, we suppose that  $$\text{ $d_{G_0}(w_1)\ge d_{G_0}(w_2)\cdots \ge d_{G_0}(w_{s-1})$  and  $d_{G_0}(z_1)\ge d_{G_0}(z_2)\cdots \ge d_{G_0}(z_{t-1})$.}$$

 \par\noindent{\bf Case 1.     $\mathcal {R}(G)=B_1$.} Since $G_0$ contains six  $\left(3^{\ge},3^{\ge }\right)$ edges and $6\left(\frac{9}{4}\right)^3+24>2\left(\frac{12}{5}\right)^3+64$,   we have   $n-7\le p\le n-6$ by \eqref{51e}.

\par\noindent{\bf Subcase 1.1. $p=n-7$.} By Lemma  \ref{52l}$(1)$,  we have  $s=t=3$.  Since $6\left(\frac{9}{4}\right)^3+16>3\left(\frac{9}{4}\right)^3+48$, we have $k\ge 5$ by \eqref{53e}. By Lemma \ref{52l}$(2)$, we can conclude that $d_{G_0}(w_2)=d_{G_0}(z_2)=3$. Suppose that $d_{G_0}(w_1)\ge \,d_{G_0}(z_1)$. By  Lemma \ref{52l}$(3)$, we have $d_{G_0}(z_1)=3$ and thus exactly one vertex  of $\{w_2,\,z_1,\,z_2\}$ is adjacent with one non-pendant vertex, as $k\ge 5$ and $p=n-7$. Combining this with     $k\ge 5$ and $p=n-7$, we have $d_{G_0}(w_1)\ge 4$ and thus    $G_0$ contains two  $\left(4^{\ge},3^{\ge }\right)$ edges and five   $\left(3^{\ge},3^{\ge }\right)$ edges. Since $5\left(\frac{9}{4}\right)^3+2\left(\frac{12}{5}\right)^3+8>2\left(\frac{12}{5}\right)^3+64$, the result follows from  \eqref{51e}.

\par\noindent{\bf Subcase 1.2. $p=n-6$.} Then, $x_0y_0\in E(G_0)$ and $s=t=3$. Since $n\ge 15$ and $5\left(\frac{9}{4}\right)^3+2\left(\frac{12}{5}\right)^3>3\left(\frac{9}{4}\right)^3+48$, we have $k\ge 5$ by \eqref{53e}.
 Combining with $k\ge 5$ and $G_0$ being a $k$-apex tree, we may suppose that $d_{G_0}(w_1)\ge d_{G_0}(w_2)\ge 4$. This implies that $G_0$ contains one $\left(4^{\ge},4^{\ge }\right)$ edge, two $\left(4^{\ge},3^{\ge }\right)$ edges and four  $\left(3^{\ge},3^{\ge }\right)$ edges, contrary with Lemma \ref{53l}$(1)$.

\par\noindent{\bf Case 2.    $\mathcal {R}(G_0)=B_2$.} Then, $G_0$ contains four $\left(4^{\ge},3^{\ge }\right)$ edges and two $\left(3^{\ge},3^{\ge }\right)$ edges. By Lemma \ref{52l}$(3)$,  we have $n-6\le p\le n-5$.

\par\noindent {\bf Subcase 2.1. $p=n-6$}. By Lemma \ref{53l} $(2)$ and $(5)$, we have $s=t=3$ and $d_{G_0}(w)=3$ holds for each vertex $w\in V\big(\mathcal {R}(G_0)\big)\backslash \{x_0\}$. Furthermore, Lemma \ref{53l}(4) implies that $d_{G_0}(x_0)=4$.  Since $n\ge 15$, we may suppose that $w_1$ is adjacent with some branching vertex of $T_{w_1}$, contrary with Lemma \ref{53l}$(5)$.
\par\noindent {\bf Subcase 2.2. $p=n-5$}. Then, $s=t=3$.    If $d_{G_0}(x_0)\ge 5$, then  $k\ge 5$ by \eqref{53e}, as $2\left(\frac{9}{4}\right)^3+4\left(\frac{15}{6}\right)^3>3\left(\frac{9}{4}\right)^3+48$. Since  $G_0$ is a $k$-apex tree with $k\ge 5$ and $p=n-5$, we may suppose that  $d_{G_0}(w_1)\ge 4$, contrary with Lemma \ref{54l}$(1)$. Otherwise,  $d_{G_0}(x_0)=4$.

   Since $n\ge 15$, we have  $d_{G_0}(w_1)\ge 5$. By Lemma \ref{54l}$(2)$, we have $d_{G_0}(w_2)=3$ and thus $d_{G_0}(w_1)=5$ by Lemma \ref{54l}$(3)$. Since $n\ge 15$, we have $d_{G_0}(z_1)\ge 5$, which is contrary with Lemma \ref{54l}$(2)$.
\end{proof}

\par\noindent{\bf Proof of Proposition \ref{51p}:}  By contradiction, assume that $G_0$ is a bicyclic graph. By Lemmas \ref{57l} and  \ref{58l}, we have   $\mathcal {R}(G_0)=Q_{s_1,s_2,s_3}$ and each vertex of every cycle of  $\mathcal {R}(G_0)$ is a branching vertex. Since $G$ is a $k$-apex tree with $k\ge 4$, we may suppose that   $d_{\mathcal {R}(G)}(x_0)=d_{\mathcal {R}(G_0)}(y_0)=3$ and $d_{G_0}(x_0)\ge d_{G_0}(y_0)\ge 4$. With the same reason, we have  $\min\big\{|V(T_{x_0})|,\,|V(T_{y_0})\big\}\ge k\ge 4$, and thus  for $w\in \big\{x_0,\,y_0\big\}$,  either $d_{G_0}(w)\ge 5$ or $d_{G_0}(w)=4$, where $d_{G_0}(w)=4$ implies that the unique neighbor vertex of $w$ in $T_w$ is a non-pendant vertex.

Let  $s=s_1+s_2-2$ and  $V(C_s)=V(P_{s_1})\cup V(P_{s_2})$.   If $p=n-4$, then $d_{G_0}(x_0)\ge d_{G_0}(y_0)\ge  6$ and thus the result follows from \eqref{51e}, as $2\left(\frac{18}{7}\right)^3+\left(\frac{36}{10}\right)^3+16>2\left(\frac{12}{5}\right)^3+64$. Thus,  $p\le n-5$.

 \par\noindent
 {\bf Case 1.  $u_0v_0\in E(G)$.} In this case, $G_0$ contains one  $\left(4^{\ge},4^{\ge }\right)$ edge and four  $\left(4^{\ge},3^{\ge }\right)$ edges. By Lemma \ref{53l}$(2)$, we have $p=n-5$.

  Since $\min\big\{|V(T_{x_0})|,\,|V(T_{y_0})|\big\}\ge k\ge 4$,  $d_{G_0}(x_0)\ge 6$ and thus $G_0$   contains one  $\left(6^{\ge},4^{\ge }\right)$ edge, two   $\left(6^{\ge},3^{\ge }\right)$ edges and two    $\left(4^{\ge},3^{\ge }\right)$ edges, contrary with Lemma \ref{54l}(3).

 \par\noindent
 {\bf Case 2.  $u_0v_0\not\in E(G)$.} In this case, $G_0$ contains six  $\left(4^{\ge},3^{\ge }\right)$ edges. If $p\le n-6$, then the result follows from Lemma \ref{21l} and \eqref{51e}, as $6\left(\frac{12}{5}\right)^3+8+3>2\left(\frac{12}{5}\right)^3+64$. Thus,     $p=n-5$.

   Since  $\min\big\{|V(T_{x_0})|,\,|V(T_{y_0})|\big\}\ge k\ge 4$, we have  $d_G(u_0)\ge d_G(v_0)\ge 6$ and  so $G_0$ contains six $\left(6^{\ge},3^{\ge }\right)$ edges. Since $6\left(\frac{18}{7}\right)^3>2\left(\frac{12}{5}\right)^3+64$, the result follows from \eqref{51e}.   \qed

\section{ $G_0$ is not a tricyclic graph}
In what follows, we shall show that $G_0$ is not a tricyclic. Throughout this section, we assume that  $G_0$ is  a tricyclic graph.
By Lemma \ref{45l}, we have $$\text{$n-8\le p\le n-4$ and $G_0$ contains a special cycle $C_s$}.$$
\begin{proposition}\label{61p} Let $G_0$ be a graph in $\mathbb{T}^k_n$ minimizing the AZI. Then $G_0$ is not tricyclic.\end{proposition}

\begin{lemma}\label{61l} $\mathcal{R}(G_0)\ne K_4$ and  $p\le n-5$.
\end{lemma}
\begin{proof} Since $G_0$ is a tricyclic graph, it suffices to show that $\mathcal{R}(G_0)\ne K_4$.  Assume that      $\mathcal {R}(G_0)=K_4$.  Since $G_0$ is a $k$-apex tree with $k\ge 4$,  $K_4$ contains at least three vertices with degree at least four in $G_0$. Thus, $G_0$ contains   three   $\left(4^{\ge},3^{\ge }\right)$-edges and  three   $\left(4^{\ge},4^{\ge }\right)$-edges. By Corollary \ref{23l}, we have
\begin{align*}AZI(G_0)\ge 3\left(\frac{12}{5}\right)^3+3\left(\frac{16}{6}\right)^3+(n-4)\left(\frac{n-1}{n-2}\right)^3. \end{align*}
Since $3\left(\frac{12}{5}\right)^3+3\left(\frac{16}{6}\right)^3>2\left(\frac{12}{5}\right)^3+64$, the result follows from \eqref{51e}. \end{proof}

\begin{lemma}\label{62l} If $k=4$, then   $p=n-5$ when one of the followings holds:\par\noindent
 $(1)$ $G_0$ contains two  $\left(5^{\ge},3^{\ge }\right)$ edges and one $\left(3^{\ge},3^{\ge }\right)$ edge;\par\noindent
 $(2)$ $G_0$  contains two    $\left(4^{\ge},3^{\ge }\right)$ edges and two $\left(3^{\ge},3^{\ge }\right)$ edges; \par\noindent
 $(3)$ $G_0$  contains five    $\left(3^{\ge},3^{\ge }\right)$ edges or three  $\left(4^{\ge},3^{\ge }\right)$ edges or two $\left(6^{\ge},3^{\ge }\right)$ edges.\end{lemma}
\begin{proof} By Lemma \ref{61l}, we assume that $p\le n-6$.  Since $2\left(\frac{15}{6}\right)^3+ \left(\frac{9}{4}\right)^3+40+2>2\left(\frac{18}{7}\right)^3+48+2>3\left(\frac{12}{5}\right)^3+40+2>5\left(\frac{9}{4}\right)^3+24+2>2\left(\frac{12}{5}\right)^3+ 2\left(\frac{9}{4}\right)^3+32>3\left(\frac{9}{4}\right)^3+48$, contrary with  \eqref{53e}.       \end{proof}

Recall that $G^*_0$ is the graph obtained from $G_0$ by deleting all the pendant vertices of $G_0$.
\begin{lemma}\label{65l}  Let $G_0$ be a tricyclic graph with $u_0\in \mathcal {R}(G_0)$. \par\noindent
   \par\noindent   $(1)$ Let $u_0$ be  a vertex of  $G^*_0$. If $G^*_0-u_0$ is a tree, then   $d_{G_0}(u_0)\ge d_{G^*_0}(u_0)+k-1\ge d_{G^*_0}(u_0)+3$;  \par\noindent
$(2)$   Let    $u_0v_0$ be a $\left(3^{\ge},3^{\ge }\right)$ edge of  $G^*_0$ such that $u_0x_0\in E(G^*_0)$ and $v_0y_0\in E(G^*_0)$, where $x_0$ and $y_0$ are both branching vertices of $G^*_0$ and $x_0=y_0$ is also permitted. If $p=n-6$ and $G^*_0-\{u_0,v_0\}$ is a tree,  then  $\min\big\{d_{G^*_0}(x_0), d_{G^*_0}(y_0)\big\}=3$ and $x_0y_0\not\in E(G^*_0)$.    \end{lemma}
\begin{proof}  Since $G_0$ is a $k$-apex tree, if  $G_0-V(T_w)$ is a tree for some $w\in \mathcal {R}(G_0)$, then $|V(T_w)|\ge k\ge 4$. Thus,   $(1)$ hold clearly. Next we show  $(2)$ holds.    Assume that $\min\big\{d_{G^*_0}(x_0), d_{G^*_0}(y_0)\big\}\ge 4$. If $d_{G_0}(v_0)\ge d_{G_0}(u_0)\ge 4$, then $G_0$ contains three  $\left(4^{\ge},4^{\ge }\right)$ edges. Since $3\left(\frac{16}{6}\right)^3+43>2\left(\frac{12}{5}\right)^3+64$, we can get a contradiction with   \eqref{51e}. Thus, we may suppose that $d_{G_0}(u_0)=3$ and thus $d_{G_0}(v_0)\ge 5$, as $G_0$ is a $k$-apex tree with $k\ge 4$. By Lemma \ref{62l}$(1)$,  we have  $k\ge 5$ and thus $d_{G_0}(u_0)=3$ and  $d_{G_0}(v_0)\ge 6$. This implies that $G_0$ contains one  $\left(6^{\ge},4^{\ge }\right)$ edge  and one  $\left(6^{\ge},3^{\ge }\right)$ edge. Since $\left(\frac{18}{7}\right)^3+\left(\frac{24}{8}\right)^3+48>2\left(\frac{12}{5}\right)^3+64$, contrary with  \eqref{51e}.  Now, we can conclude that $\min\big\{d_{G^*_0}(x_0), d_{G^*_0}(y_0)\big\}=3$.

 If $x_0y_0\in E(G^*_0)$, then either $\min\big\{d_{G_0}(u_0),d_{G_0}(v_0)\big\}\ge 4$ or   $\min\big\{d_{G_0}(u_0),d_{G_0}(v_0)\big\}=3$.
When $\min\big\{d_{G_0}(u_0),d_{G_0}(v_0)\big\}\ge 4$, then $G_0$ contains one  $\left(4^{\ge},4^{\ge }\right)$ edge, two $\left(4^{\ge},3^{\ge }\right)$ edges and  one  $\left(3^{\ge},3^{\ge }\right)$ edge. Since $2\left(\frac{12}{5}\right)^3+\left(\frac{9}{4}\right)^3+\left(\frac{16}{6}\right)^3+35>2\left(\frac{12}{5}\right)^3+64$, we can get a contradiction with  \eqref{51e}. When $\min\big\{d_{G_0}(u_0),d_{G_0}(v_0)\big\}=3$, then   $\max\big\{d_{G_0}(u_0),d_{G_0}(v_0)\big\}\ge 5$, as $G_0$ is a $k$-apex tree with $k\ge 4$. By Lemma \ref{62l}(1), we have  $k\ge 5$  and thus $\min\big\{d_{G_0}(u_0),d_{G_0}(v_0)\big\}=3$ and   $\max\big\{d_{G_0}(u_0),d_{G_0}(v_0)\big\}\ge 6$. This implies that  $G_0$ contains two  $\left(6^{\ge},3^{\ge }\right)$ edges  and two $\left(3^{\ge},3^{\ge }\right)$ edges. Since $2\left(\frac{18}{7}\right)^3+2\left(\frac{9}{4}\right)^3+35>2\left(\frac{12}{5}\right)^3+64$ and \eqref{51e}, we have  $x_0y_0\not\in E(G^*_0)$.  \end{proof}
\begin{figure}[h!]
\begin{center} \includegraphics[scale=0.70]{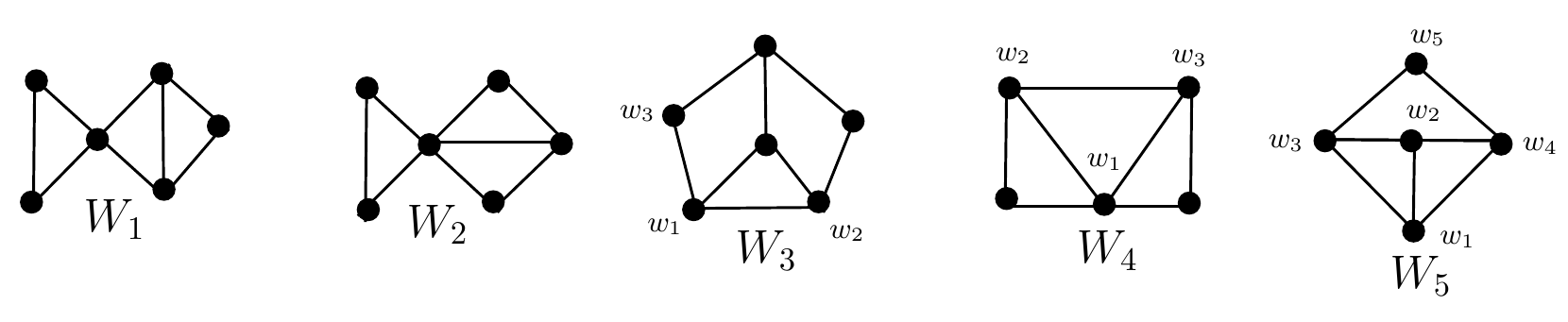}
    \caption{The  tricyclic   graphs $W_1$, $W_2$,  $\ldots$, $W_5$.}  \label{22f}\end{center}\end{figure}
    In what follows, let $W_1$, $W_2$, $\ldots$, $W_5$ be the tricyclic  graphs as shown in Fig. \ref{22f}.
\begin{lemma}\label{66l}  Any cycle of    $G_0$ contains no  $2$-vertex  of $G_0$ and $n-7\le p\le n-5$.   \end{lemma}
\begin{proof} If any  cycle of    $G_0$ contains no  $2$-vertex  of $G_0$, then   $G_0$ contains four  $\left(3^{\ge},3^{\ge }\right)$ edges. Since $4\left(\frac{9}{4}\right)^3+48>2\left(\frac{12}{5}\right)^3+64$, we have $p\ge n-7$ by \eqref{51e}. Combining this with Lemma \ref{61l}, we have $n-7\le p\le n-5$. Thus, to complete the proof of this result, it suffices to show that any cycle of    $G_0$ contains no  $2$-vertex  of $G_0$.

 By contradiction, assume that $u_0v_0$ is an edge of some cycle of $G_0$ such that $u_0$ is a  2-vertex of $G_0$.   Similar with  Lemma \ref{56l}, we have $k=4$ by Proposition \ref{51p}. Combining $k=4$ with    Lemmas \ref{23l} and  \ref{24l}, we may suppose that  $G_0-u_0v_0\in \mathscr{U}^3_n$ by  the choice of $G_0$.

 Since $C_s$ is a special cycle of $G_0$,    $G_0$  contains three $\left(3^{\ge},3^{\ge }\right)$ edges.  If $p= n-7$,  then  Corollary \ref{21c} implies that     \begin{align*} AZI(G_0)\ge 3\left(\frac{9}{4}\right)^3+48+(n-7)\left(\frac{n-1}{n-2}\right)^3,\end{align*} and thus the result follows from \eqref{53e}.  In what follows, we assume that  $$n-6\le p\le n-5.$$

\par\noindent {\bf Case 1. $\mathcal {R}(G_0)=Q_{s_1,s_2,s_3,s_4}$.} Let $x_0,y_0$ be the two 4-vertices of $Q_{s_1,s_2,s_3,s_4}$ such that   $d_{G_0}(x_0)\ge d_{G_0}(y_0)$. Since $n-6\le p\le n-5$, we have $d_{G_0}(x_0)\ge 7$. Thus,   $G_0$ contains two     $\left(7^{\ge},3^{\ge }\right)$ edges and one $\left(4^{\ge},3^{\ge }\right)$ edge. Since $\left(\frac{12}{5}\right)^3+ 2\left(\frac{21}{8}\right)^3+32+2>3\left(\frac{9}{4}\right)^3+48$,  the result follows from \eqref{53e}.

 \par\noindent {\bf Case 2.    $\mathcal {R}(G_0)$ contains a cut vertex.} Since $\mathcal {R}(G_0)$ is a tricyclic graph and $n-6\le p\le n-5$, we have  $|V(\mathcal {R}(G_0))|=6$. In this case, $\mathcal {R}(G_0)\in \big\{W_1,\,W_2\big\}$ and $p=n-6$.

If $\mathcal {R}(G_0)=W_1$, then the two 3-vertices of $W_1$ are also two 3-vertices of $G_0$ by Lemma \ref{62l}$(3)$. Combining this with $G_0$ being  a 4-apex tree, we can see that  at least one of three 2-vertices of $W_1$ is a branching vertex of $G_0$. Thus,  $G_0$  contains two     $\left(4^{\ge},3^{\ge }\right)$ edges and two $\left(3^{\ge},3^{\ge }\right)$ edges, contrary with Lemma \ref{62l}(2).

If $\mathcal {R}(G_0)=W_2$, then $G_0$ contains two     $\left(5^{\ge},3^{\ge }\right)$ edges and one $\left(3^{\ge},3^{\ge }\right)$ edge, as $C_s$ is a special cycle of $G_0$. This is contrary with Lemma \ref{62l}(1).

 \par\noindent {\bf Case  3.   $\mathcal {R}(G_0)$ contains no  cut vertex.}  By Case 1 and  Lemma \ref{61l}, we have  $$\text{ $\mathcal {R}(G_0)\not\in \big\{K_4,\,Q_{s_1,s_2,s_3,s_4}\big\}$.}$$
 \par\noindent {\bf Subcase 3.1. $p=n-6$.}  Let $u_0$ be a vertex of  $G^*_0$    with $d_{G^*_0}(u_0)\ge 4$. If   $G^*_0-u_0$ is a tree, then $d_{G_0}(u_0)\ge 7$  by Lemma \ref{65l}$(1)$. By Lemma \ref{62l}$(3)$,   $u_0$ has at most  one branching vertex as its neighbor    in $G_0$, and thus it has at least three   2-vertices of $G_0$ as its neighbors in $G^*_0$. This implies that  we can obtain  a tree from $G_0$ be deleting three neighbors of $u_0$ in $G^*_0$, a contradiction.

 Thus, for any vertex $u$ of $G^*_0$, either $d_{G^*_0}(u)\le 3$ or $G^*_0-u$ is not a tree.   Combining  this  with  Lemma \ref{65l}$(2)$, we can conclude that $G^*_0=W_3$ by the table of connected graphs with   six vertices (see  \cite{Dc1}).  By Lemma \ref{62l}(2), we have $d_{G_0}(w_1)=d_{G_0}(w_2)=3$ and thus   $G_0-\big\{w_1,w_2\big\}$ is a tree, a contradiction.

 \par\noindent {\bf Subcase 3.2. $p=n-5$.} Since $\mathcal {R}(G_0)\not\in \big\{K_4,\,Q_{s_1,s_2,s_3,s_4}\big\}$, we have $G^*_0\in \big\{W_4,\,W_5\big\}$ by the table of connected graph with five vertices (see    \cite{Dc1}).

If  $G^*_0=W_4$, then $d_{G_0}(w_1)\ge 7$ by Lemma \ref{65l}$(1)$. Since $G_0$ is a 4-apex tree, we may suppose that $d_{G_0}(w_2)\ge 4$. In this case, $G_0$ contains one  $\left(7^{\ge},4^{\ge }\right)$ edge  and  one  $\left(7^{\ge},3^{\ge }\right)$ edge. Since $\left(\frac{28}{9}\right)^3+\left(\frac{21}{8}\right)^3+40>3\left(\frac{9}{4}\right)^3+48$, the result follows from \eqref{53e}.

  If $G^*_0=W_5$, then either  $d_{G_0}(w_1)\ge d_{G_0}(w_2)\ge 4$ or   $d_{G_0}(w_1)\ge 5>d_{G_0}(w_2)=3$, as $G_0$ is a 4-apex tree. In both cases,  $G_0$ contains
three  $\left(4^{\ge},3^{\ge }\right)$ edges and two $\left(3^{\ge},3^{\ge }\right)$ edges. Since $3\left(\frac{12}{5}\right)^3+2\left(\frac{9}{4}\right)^3+16+3>3\left(\frac{9}{4}\right)^3+48$, the result follows from \eqref{53e}.  \end{proof}

\begin{lemma}\label{67l}     $\mathcal {R}(G_0)$ contains no  cut vertex and $\mathcal {R}(G_0)\ne Q_{s_1,s_2,s_3,s_4}$. \end{lemma}
\begin{proof} By    Lemma \ref{66l},    any cycle of $G_0$ contains no 2-vertex of $G_0$. By contradiction, we assume that either  $\mathcal {R}(G_0)$ contains a  cut vertex or  $\mathcal {R}(G_0)=Q_{s_1,s_2,s_3,s_4}$.

\par\noindent {\bf Case 1.    $\mathcal {R}(G_0)$ contains a  cut vertex.} Since $\mathcal {R}(G_0)$ contains a  cut vertex and $\mathcal {R}(G_0)$ is a tricyclic graph, we have    $|V(\mathcal {R}(G_0))|\ge 6$ and thus    $n-7\le p\le n-6$ by Lemma \ref{66l}.

Since $\mathcal {R}(G_0)$ contains a  cut vertex and $\mathcal {R}(G_0)$  is a tricyclic graph, $\mathcal {R}(G_0)$  contains a cycle $C_q$ such that each vertex of  $V(C_q)\setminus \{x_0\}$ is a 2-vertex of $\mathcal {R}(G_0)$, where $x_0$ is a cut vertex of $\mathcal {R}(G_0)$. This implies that $\mathcal {R}(G_0)-\big(V(C_q)\setminus \{x_0\}\big)$ is a bicyclic graph containing at least four vertices and thus $G_0$ contains eight  $\left(3^{\ge},3^{\ge }\right)$ edges. Since $8\left(\frac{9}{4}\right)^3+2
 >2\left(\frac{12}{5}\right)^3+64$,   the result follows from \eqref{51e}.

\par\noindent {\bf Case 2.  $\mathcal {R}(G_0)=Q_{s_1,s_2,s_3,s_4}$.} Suppose that  $d_{\mathcal {R}(G_0)}(x_0)=d_{\mathcal {R}(G_0)}(y_0)=4$. By Lemma  \ref{66l},  $n-7\le p\le n-5$.  Since $G_0$ is a $k$-apex tree with $k\ge 4$, we have $\min\big\{|V(T_{x_0})|, |V(T_{y_0})|\big\}\ge k\ge 4$, and thus $d_{G_0}(x_0)\ge d_G(y_0)\ge 5$. Since no cycle of $G_0$ contains a 2-vertex of $G_0$,  $G_0$ contains seven   $\left(5^{\ge},3^{\ge }\right)$ edges. Now, the result follows from \eqref{51e}, as  $7\left(\frac{15}{6}\right)^3>2\left(\frac{12}{5}\right)^3+64$.   \end{proof}

\par\noindent{\bf Proof of Proposition \ref{61p}:} By contradiction, assume that   $G_0$ is  a tricyclic graph.  From  Lemmas  \ref{61l}, \ref{66l}  and   \ref{67l}, it follows that   \begin{align}\label{61e}\text{$n-7\le p\le n-5$, $\mathcal {R}(G_0)\not\in \big\{Q_{s_1,s_2,s_3,s_4},\,\,K_4\big\}$ and $\mathcal {R}(G_0)$ contains no cut vertex}.\end{align}

 If $6\le |V(\mathcal {R}(G_0))|\le 7$, then $\mathcal {R}(G_0)$ contains at least eight  $\left(3^{\ge},3^{\ge }\right)$ edges, as $\mathcal {R}(G_0)$ contains no cut vertex and each vertex of $\mathcal {R}(G_0)$ is a branching vertex by Lemmas \ref{66l} and \ref{67l}. Since $8\left(\frac{9}{4}\right)^3+2>2\left(\frac{12}{5}\right)^3+64$, the result follows from \eqref{51e}. Next, suppose that $|V(\mathcal {R}(G_0))|=5$, as $\mathcal {R}(G_0)\ne K_4$. Combining $|V(\mathcal {R}(G_0))|=5$ and $\mathcal {R}(G_0)\ne Q_{s_1,s_2,s_3,s_4}$, we have $\mathcal {R}(G_0)\in \big\{W_4,\, W_5\big\}$ by the table of connected graphs with five  vertices (see   \cite{Dc1}). In this case, $G_0$ contains at least seven $\left(3^{\ge},3^{\ge }\right)$ edges. Since $7\left(\frac{9}{4}\right)^3+16> 2\left(\frac{12}{5}\right)^3+64$ and \eqref{51e}, we have $n-6\le p\le n-5$.

If   $\mathcal {R}(G_0)=W_4$, then $\mathcal {R}(G_0)$ contains four $\left(4^{\ge},3^{\ge }\right)$ edges and three $\left(3^{\ge},3^{\ge }\right)$ edges, as   $W_4$ contains no 2-vertex of $G_0$.    Since $4\left(\frac{12}{5}\right)^3+3\left(\frac{9}{4}\right)^3+3> 2\left(\frac{12}{5}\right)^3+64$ and $p\ge n-6$, the result follows from \eqref{51e}.

  Otherwise,   $\mathcal {R}(G_0)=W_5$.    Since $G_0$ is $k$-apex tree with $k\ge 4$, we have $\max\{d_{G_0}(w_1),\,d_{G_0}(w_2)\}$\\$\ge 4$, which implies that $G_0$ contains three  $\left(4^{\ge},3^{\ge }\right)$ edges and four $\left(3^{\ge},3^{\ge }\right)$ edges.   Since $3\left(\frac{12}{5}\right)^3+4\left(\frac{9}{4}\right)^3+8>
2\left(\frac{12}{5}\right)^3+64$ and \eqref{51e}, we have $p=n-5$.

   If either $d_{G_0}(w_3)=3$ or $d_{G_0}(w_4)=3$, then $d_{G_0}(w_1)\ge d_{G_0}(w_2)\ge 5$, as $G_0$ is a $k$-apex tree with $k\ge 4$. Thus,  $G_0$ contains
five  $\left(5^{\ge},3^{\ge }\right)$ edges.  Otherwise,   $d_{G_0}(w_3)\ge d_{G_0}(w_4)\ge 4,$ which implies that   $G_0$ contains six $\left(4^{\ge},3^{\ge }\right)$ edges and one $\left(3^{\ge},3^{\ge }\right)$ edge.

Since $6\left(\frac{12}{5}\right)^3+\left(\frac{9}{4}\right)^3> 5\left(\frac{15}{6}\right)^3+16>
2\left(\frac{12}{5}\right)^3+64$, the result follows from \eqref{51e}. \qed

\par\medskip
\noindent
 {\it Acknowledgement\/}.   The first  author is supported by the  Natural Science Foundation of Guangdong Province (No. 2022A1515011786).

\end{document}